\documentclass[a4paper,10pt,final]{elsarticle}
\usepackage[utf8]{inputenc}
\usepackage{import}

\usepackage{fullpage}

\usepackage[english]{babel}

\usepackage{amsmath,amssymb,amsthm}
\usepackage{cases}
\usepackage{mathtools}
\usepackage{xfrac}

\usepackage{siunitx}
\usepackage{pdflscape}
\usepackage{tabularx}
\usepackage{booktabs}         

\usepackage{multicol}
\usepackage{multirow}
\usepackage{xparse}
\usepackage{xspace}
\usepackage{makecell}

\usepackage[usenames,dvipsnames,table]{xcolor}

\usepackage{footnote} 
\usepackage{graphicx}
\usepackage{placeins}
\usepackage{overpic}
\usepackage{svg}
\svgpath{{./images/}}

\definecolor{ultramarine}{RGB}{0,32,96}
\usepackage[colorlinks=true]{hyperref}         
\usepackage{cleveref}
\usepackage{caption}          
\usepackage[notref]{showkeys}
\usepackage[color=gray, backgroundcolor=gray!10, textwidth=2cm, textsize=footnotesize]{todonotes}
\usepackage{setspace}
\usepackage{subcaption}

\usepackage{framed} 

\usepackage{enumitem}

\newcommand{\Sec}[1]{Sec.~\ref{#1}}
\newcommand{\Fig}[1]{Fig.~\ref{#1}}
\newcommand{\Tab}[1]{Tab.~\ref{#1}}

\newcommand{\bra}[1]{\left( #1 \right)}

\newcommand{\sbra}[1]{\left[ #1 \right]}
\newcommand{\cbra}[1]{\left\lbrace #1 \right\rbrace}

\renewcommand{\l}[1]{\text{\scalebox{0.9}{\upshape #1}}}

\makeatletter
\newcommand{\dotr}[1]{%
  \mathpalette\@dotr{#1}%
}
\newcommand*{\@dotr}[2]{%
  \sbox0{$\m@th#1#2$}%
  \usebox{0}%
  \raisebox{\dimexpr\ht0-\height}{$\m@th#1\@smallbullet#1\bullet$}%
  \kern\scriptspace
}
\newcommand*{\@smallbullet}[2]{%
  \scalebox{.5}{$\m@th#1#2$}%
}
\makeatother

\newcommand{\varepsilonb}{\boldsymbol{\varepsilon}}

\newcommand{\Dm}{\mathcal{D}}
\newcommand{\Em}{\mathcal{E}}

\newcommand{\Pm}{\mathcal{P}}

\newcommand{\gsf}{\mathsf{g}}
\newcommand{\hsf}{\mathsf{h}}
\newcommand{\lsf}{\mathsf{l}}

\newcommand{\wsf}{\mathsf{w}}

\newcommand{\idest}{\emph{i.e.}\xspace}

\newcommand{\coh}{\frac{1}{2}}
\newcommand{\ar}[1]{\!\bra{#1}}
\newcommand{\abs}[1]{\lvert #1 \rvert}
\newcommand{\diff}{\mathop{}\!\mathrm{d}}
\newcommand{\dx}{\diff\spos}
\renewcommand{\t}{t}
\NewDocumentCommand\disspot{o}{\IfNoValueTF{#1}{\Dm}{\varphi\ar{#1}}}

\NewDocumentCommand\Enint{o}{\IfNoValueTF{#1}{\Em}{\Em\ar{#1}}} 
\NewDocumentCommand\Enpot{o}{\IfNoValueTF{#1}{\Pm}{\Pm\ar{#1}}}
\NewDocumentCommand\intargt{O{0}O{\t}m}{\int_{#1}^{#2} #3 \,\diff \tau}
\newcommand{\spos}{x}
\NewDocumentCommand\sstr{o}{\IfNoValueTF{#1}{\sigma}{\sigma\ar{#1}}}
\NewDocumentCommand\sstra{o}{\IfNoValueTF{#1}{\tilde\sigma}{\tilde\sigma\ar{#1}}}
\NewDocumentCommand\teps{o}{\IfNoValueTF{#1}{\varepsilonb}{\varepsilonb\ar{#1}}}
\NewDocumentCommand\sdamdt{o}{\IfNoValueTF{#1}{\dot{\sdam}}{\dot{\sdam}\ar{#1}}}
\newcommand{\fDeg}{\gsf_{_{\scalebox{0.7}{$\sintL$}}}} 
\newcommand{\fCom}{\hsf_{_{\scalebox{0.7}{$\sintL$}}}} 
\newcommand{\fDam}{\wsf} 
\newcommand{\lDam}{\lsf} 
\newcommand{\sdisp}{u} 
\NewDocumentCommand\seps{o}{\IfNoValueTF{#1}{\epsilon}{\epsilon\ar{#1}}}
\NewDocumentCommand\sdam{o}{\IfNoValueTF{#1}{\alpha}{\alpha\ar{#1}}}
\NewDocumentCommand\sdammax{o}{\IfNoValueTF{#1}{\bar\alpha}{\bar\alpha\ar{#1}}}
\newcommand{\sope}{\delta} 
\newcommand{\sintL}{\ell} 
\newcommand{\scristr}{\bar{\sstr}}
\newcommand{\scriope}{\bar{\sope}}
\newcommand{\scrieps}{\bar{\seps}}
\newcommand{\cw}{c_{\fDam}} 

\newcommand{\enfree}{\psi}
\newcommand{\enfrac}{\gamma}

\newcommand{\Gc}{G_{\l{c}}} 
\newcommand{\Gcl}{G_{\l{c\,l}}} 
\newcommand{\GcII}{G_{\l{c\,2}}} 

\definecolor{DarkerGreen}{RGB}{0,170,0}
\definecolor{DarkerRed}{RGB}{170,0,0}

\definecolor{myRed}{rgb}{0.450385, 0.157961, 0.217975}
\definecolor{myBlue}{rgb}{0.139681, 0.311666, 0.550652}

\definecolor{myKcolor}{rgb}{0, 0.411765, 0.572549}
\definecolor{myEcolor}{rgb}{0, 0, 0}

\newdefinition{rem}{Remark}

\newcommand{\damBar}[1][$\sdammax$]{
\colorbox{white}{
\begin{minipage}{12mm}
  \fontsize{8pt}{9.6pt}\selectfont 
  \centering
  #1\\
  $0$\:\raisebox{0.2mm}{\includegraphics[width=8mm]{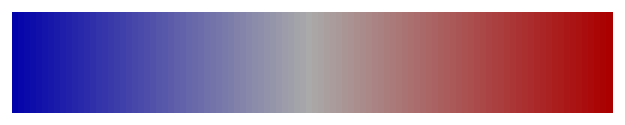}}\:$1$
\end{minipage}
}}

\newcommand{\LLeg}{
\colorbox{white}{
\begin{minipage}{12mm}
\fontsize{7pt}{8.4pt}\selectfont 
\raggedright
  \includegraphics[trim=2mm 1.8mm 0mm 0mm,clip,scale=0.8]{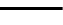}\:\eqref{mod_L1} \\
  \includegraphics[trim=2mm 1.8mm 0mm 0mm,clip,scale=0.8]{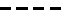}\:\eqref{mod_L2} \\
  \:\includegraphics[trim=2mm 1.8mm 0mm 0mm,clip,scale=0.8]{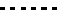}\eqref{mod_L12}
\end{minipage}
}}

\newcommand{\lLeg}[1][1]{
\colorbox{white}{
\begin{minipage}{20mm}
\fontsize{7pt}{8.4pt}\selectfont 
  $(\sdammax = #1)$ \\[1ex]
  \includegraphics[trim=2mm 1.8mm 0mm 0mm,clip,scale=0.8]{L1}\:$\sintL$ = \SI{10}{mm} \\
  \includegraphics[trim=2mm 1.8mm 0mm 0mm,clip,scale=0.8]{L2}\:$\sintL$ = \SI{5}{mm} \\
  \includegraphics[trim=2mm 1.8mm 0mm 0mm,clip,scale=0.8]{L3}\:$\sintL$ = \SI{1}{mm} 
\end{minipage}
}}

\newcommand{\llLeg}[1][1]{
\colorbox{white}{
\begin{minipage}{21mm}
\fontsize{7pt}{8.4pt}\selectfont 
  $(\sdammax = #1)$ \\[1ex]
  \includegraphics[trim=2mm 1.8mm 0mm 0mm,clip,scale=0.8]{L1}\:$\sintL$ = \SI{5}{mm} \\
  \includegraphics[trim=2mm 1.8mm 0mm 0mm,clip,scale=0.8]{L2}\:$\sintL$ = \SI{2.5}{mm} \\
  \includegraphics[trim=2mm 1.8mm 0mm 0mm,clip,scale=0.8]{L3}\:$\sintL$ = \SI{0.5}{mm} 
\end{minipage}
}}

\newcommand{\lDILeg}[1][1]{
\colorbox{white}{
\begin{minipage}{22mm}
\fontsize{7pt}{8.4pt}\selectfont 
  $(\sdammax = #1)$ \\[1ex]
  \includegraphics[trim=2mm 1.8mm 0mm 0mm,clip,scale=0.8]{L1}\:$\sintL$ = \SI{1}{mm} \\
  \includegraphics[trim=2mm 1.8mm 0mm 0mm,clip,scale=0.8]{L2}\:$\sintL$ = \SI{0.5}{mm} \\
  \includegraphics[trim=2mm 1.8mm 0mm 0mm,clip,scale=0.8]{L3}\:$\sintL$ = \SI{0.1}{mm} 
\end{minipage}
}}

\newcommand{\HLeg}{
\colorbox{white}{
\begin{minipage}{12mm}
\fontsize{7pt}{8.4pt}\selectfont 
  \includegraphics[trim=2mm 1.8mm 0mm 0mm,clip,scale=0.8]{L1}\:\eqref{mod_H1} \\
  \includegraphics[trim=2mm 1.8mm 0mm 0mm,clip,scale=0.8]{L2}\:\eqref{mod_H2}
\end{minipage}
}}

\newcommand\eps{\varepsilon}

\numberwithin{equation}{section}
\usepackage{fancyhdr}
\newcommand\Functeps{\mathcal F_\eps}

\newcommand\PotDann{\omega}
\newcommand\FailureS{\overline{\varsigma}}

\newcommand\sfratt{s_{\mathrm{frac}}}

\newcommand\hatf{l}

\newcommand\kscal{k_0}
\DeclareMathOperator*{\argmax}{arg\,max}

\graphicspath{
  {images/}
  {images/mathematica/}}

\journal{~}

\begin{document}

\begin{frontmatter}

\title{Phase-field modelling of cohesive fracture. Part III: From mathematical results to engineering applications}

\author[unipi]{Roberto Alessi}
\author[unifi]{Francesco Colasanto}
\author[unifi]{Matteo Focardi}  

\address[unipi]{DICI, Università di Pisa, Italy}
\address[unifi]{DIMAI, Università di Firenze, Italy}

\begin{abstract}
This paper concludes a three-part effort aimed at developing a consistent and unified framework for the phase-field modeling of cohesive fracture. Building on the theoretical foundations established in the first two parts, which included a $\Gamma$-convergence result for a broad class of phase-field energy functionals and the presentation of a rigorous analytical methodology for constructing models tailored to specific cohesive laws, this third paper explores the mechanical response of phase-field models, most of which are novel, associated with different cohesive fracture behaviors within a one-dimensional framework.
Particular emphasis is placed on the possibility of formulating distinct phase-field models that, despite exhibiting different evolutions of their phase-field and displacement profiles, yield identical cohesive fracture responses.
Thus, this work aims at providing a practical interpretation of the mathematical framework connecting  the theoretical insights established in the previous parts for physical relevant applications.
\end{abstract}

\end{frontmatter}

\tableofcontents

\clearpage

%

\section{Introduction}
\label{}

\subsection{Background and Motivation}
\label{sec_bac_and_mot} 
 
Predicting fractures and their propagation is of utmost importance for allowing the design of safer and more resilient mechanical systems and preventing structural failures.
Griffith's brittle fracture theory, introduced in a seminal article \cite{Griffith1921}, remains one of the most widely used and simplest fracture theories to date. 
In Griffith's model, the surface fracture energy is assumed to dissipate entirely upon the creation of a new discontinuity (the crack), irrespective of the magnitude of the displacement jump across the crack surfaces.
This implies the absence of forces between the crack surfaces, regardless of the crack opening.
Because of its simplicity, Griffith's brittle fracture model owns two major flaws such as (i) its inability to account for initiation of a crack in a pristine elastic body (without a pre-existing notch), \cite{marigo2010,Tanne2018,Kumar2020}, (ii) its prediction of unrealistic scale effects, \cite{Marigo2023}.

Cohesive fracture models address these limitations by considering non-vanishing forces, the cohesive forces, between the crack surfaces. Consequently, the surface energy density becomes a function of the displacement jump, an idea first proposed in \cite{Dugdale1960,Barenblatt1962}. Typical cohesive laws, such as those introduced by Dugdale and Barenblatt, are characterized by a critical stress and a characteristic length. These properties remedy the deficiencies of Griffith's model and offer a more realistic representation of fracture processes, \cite{Marigo2023}. 

The reformulation of Griffith's brittle fracture theory as a free-discontinuity energy minimization problem, that turned out to be a crucial step towards the development of phase-field fracture models, has been carried out in \cite{Marigo1998,Focardi2001}. A similar variational interpretation for cohesive fracture models has been accomplished in \cite{Bourdin2008}. 

The brittle fracture problem, formulated as a free-discontinuity energy minimization problem, can be regularized through a phase-field approach \cite{Bourdin2000b}, inspired by the work of \cite{Ambrosio1992}. This regularization, framed within the variational theory of $\Gamma$-convergence, allows the crack path to emerge naturally from the localization of a smooth phase-field as the regularization length tends to zero \cite{Braides1998}. 

This regularization has enabled the numerical implementation of phase-field brittle fracture models, allowing engineers to simulate complex fracture processes that were previously unfeasible with classical tools \cite{Bourdin2014,Mesgarnejad2015}. 
The phase-field approach to fracture is undoubtedly one of the leading modeling strategies in the fracture mechanics field since it allows, for instance, to capture the nucleation of cracks with or without pre-existing singularities and describe complex crack patterns, offering also a straightforward numerical implementation, typically based on alternate minimization schemes.
The contribution of \cite{Bourdin2000b} can be therefore considered as a fundamental bridging work between the mathematical community of calculus of variations and the engineering community of fracture mechanics.

Under the mechanical viewpoint, the phase-field variable can be interpreted as a damage variable~\cite{LemaitreB1985} and phase-field models as gradient damage models, \cite{Pham2011}. Indeed, the phase-field approach to brittle fracture has promoted a flourishing and still ongoing research activity on gradient-damage models, 
\cite{Pham2010c,Pham2010a,Miehe2010a,Marigo2016}.

Phase-field models for brittle fracture inherit a significant limitation from Griffith's theory: they do not allow for independent control of critical strength, fracture toughness (or critical energy release rate), and regularization length. As a result, the fracture nucleation stress threshold cannot be directly linked with Griffith's sharp-interface model. This limitation restricts the development of a flexible and general framework capable of accurately modeling crack nucleation in both smooth and notched domains~\cite{Bourdin2008,Tanne2018,Lopez-Pamies2024}.

As highlighted in \cite{Marigo2023}, cohesive fracture models appear to be the most natural candidates to overcome the limitations of Griffith's brittle fracture theory.
Inspired by \cite{Pham2010c,Pham2010a}, a standalone gradient-damage (phase-field) model, that is a model not deduced as a regularization of a free-discontinuity fracture model, was developed by \cite{Lorentz2011}, where a cohesive fracture failure response in a one-dimensional uniaxial tension setting has been investigated in a closed form. 
This model has then been extended to higher dimensional problems and numerically implemented to simulate failure processes in large scale domains \cite{Lorentz2011a}. 
Although promising, this cohesive phase-field model was initially overlooked by much of the scientific community working in the fracture mechanics field.

Meanwhile, a successful attempt to develop a mathematically consistent variational phase-field model regularizing the free-discontinuity cohesive fracture problem of \cite{Bourdin2008} was provided by \cite{Conti2016}. 
However, despite its mathematical rigor, the initial lack of numerical implementation posed challenges for its widespread adoption.
Nevertheless, this work can be considered, similarly to \cite{Bourdin2000b}, another example of fruitful interaction between the mathematical and engineering communities.
A numerical implementation of the cohesive phase-field model proposed by \cite{Conti2016} was subsequently provided by \cite{Freddi2017}. This implementation raised quite several issues, suh as the setup of a backtracking algorithm and a further regularization of the originally truncated degradation function with respect to the internal length.
Moreover, the difficulty, if not impossibility, of easily understanding how to tune the elastic degradation function to describe a specific cohesive law and the use of a fixed quadratic phase-field dissipation function, further limiting the flexibility of the mechanical response, have likely been obstacles to the adoption of this model within the fracture mechanics community. 

More impactful contributions came from \cite{Wu2017,Wu2018b}, which, in the wake of \cite{Lorentz2011}, proposed a new generation of phase-field cohesive fracture models built within the variational framework of gradient-damage models~\cite{Pham2010c,Pham2010a}.
These models incorporate polynomial crack geometric functions and rational energetic degradation functions, allowing for the independent tuning of critical stress, fracture toughness, and regularization length to match target softening laws.
The key modeling idea enabling these models \cite{Conti2016,Wu2017,Wu2018b} to describe cohesive fracture failures has been to make the elastic degradation function dependent on the regularization length.

Despite the model flexibility and soundness, how to set the material function to obtain a specific target cohesive fracture response remained rather obscure. In this regard, a crucial contribution came from \cite{Feng2021}, which introduced an integral relation linking a single unknown function defining both the degradation function and the local phase-field dissipation function to the desired traction-separation law. 
Interestingly and differently from \cite{Conti2016,Lammen2023,Lammen2025}, 
the cohesive phase-field fracture models in \cite{Wu2017,Wu2018b,Feng2021} are not affected by the numerical issues faced in \cite{Freddi2017}. Moreover, in these models, it is found that both the regularization length and mesh size have minimal influence on the global mechanical response, provided the mesh sufficiently resolves the regularization length. This new generation of cohesive phase-field fracture models has thus proven to be reliable and effective in simulating practical problems. 
Nevertheless, a mathematical proof of the convergence of the models \cite{Wu2017,Wu2018b,Feng2021} toward the free-discontinuity cohesive fracture problem and a rigorous analysis of how to identify the material functions complying with a given traction-separation law, are still missing. This three-part work (\cite{Alessi2025a}, \cite{Alessi2025b}, and the present third contribution) aims precisely to fill this gap.

It is worth mentioning that another way to describe cohesive fracture behaviors within a variational framework turned out to be the coupling of the phase-field brittle fracture model with the perfect plasticity model, as done in \cite{Alessi2014,Alessi2015,Alessi2018a}. Depending on the way damage and plasticity are coupled, it was shown in this work that different cohesive behaviors may be captured by the plastic field localizing as a Dirac measure, thus realizing the sought cohesive fracture response.

\begin{landscape}
\begin{figure}
 \centering
  {\fontsize{8pt}{9.6pt}\selectfont 
  \includeinkscape[width=\linewidth]{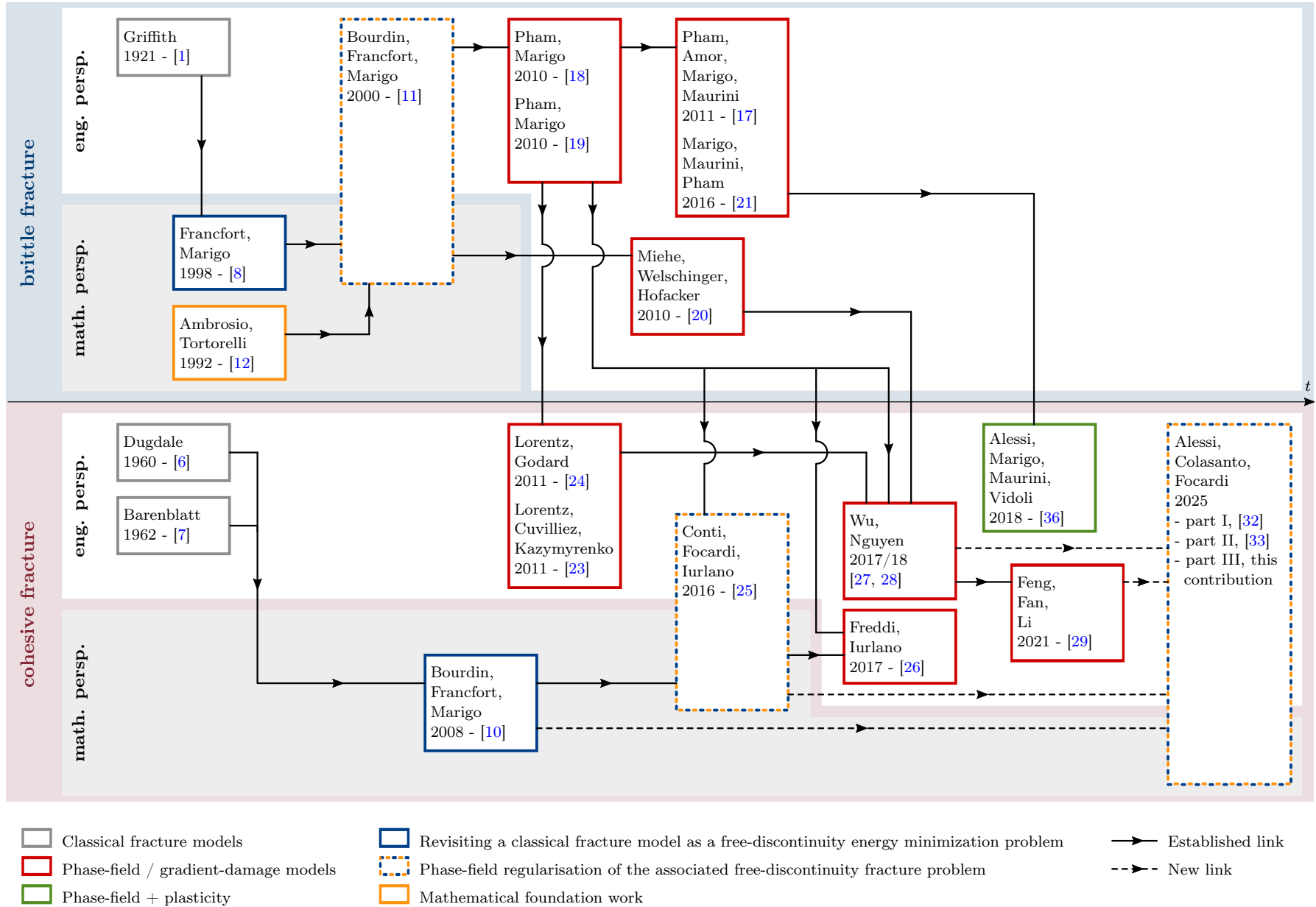}
  }
 \caption{Overview of the state of the art with both a mathematical and engineering perspective of the phase-field models describing brittle and cohesive fractures.  Specifically, milestones, their logic connection and the collocation of this three part work within the framework of the variational approach to fracture are highlighted.}
 \label{fig_FT}
\end{figure}
\end{landscape}

A visual overview of the state of the art discussed above is presented on a timeline in \Fig{fig_FT}. Therein, milestones, their logic connection and the placement of this third part contribution within the framework of the variational approach to fracture are highlighted.

\subsection{Scope of this paper}

This three-part work aims to develop a consistent and unified framework for phase-field modeling of cohesive fracture.
The previous two parts of this work, \cite{Alessi2025a,Alessi2025b}, explored mathematical aspects and provided theoretical results that lay the foundation for our final analysis. 

Specifically, the first paper establishes a $\Gamma$-convergence result for a broad class of phase-field energy functionals, encompassing models from \cite{Conti2016}, \cite{Wu2017,Wu2018b}, \cite{Feng2021} and \cite{Lammen2023,Lammen2025}. Additionally, by adjusting the functional scaling, we showed how the proposed framework also generalizes the classical Ambrosio-Tortorelli model for approximating brittle fracture energies.
The second paper develops an analytical framework for constructing phase-field models that reproduce specific cohesive laws by either fixing the degradation function and deriving the local phase-field dissipation potential, or vice versa. This approach mathematically validates and extends the results of \cite{Feng2021}, enabling the derivation of multiple phase-field models that adhere to the same target cohesive law, \idest exhibit  identical overall cohesive fracture behavior but are associated with different phase-field profile evolutions.
This third paper examines the mechanical responses of various phase-field models associated with canonical traction-separation laws in a one-dimensional setting. Adopting an engineering-oriented perspective, this study provides a link to the first two parts of this work, highlighting the theoretical and original contributions established therein.

This work is structured as follows. \Sec{sec_1Dprob} aims to bridge the mathematical setting and findings of the first two parts of this work (\cite{Alessi2025a} and \cite{Alessi2025b}) with a more engineering-oriented framework.
In this regard, key quantities and well-established results concerning the mechanical response of the one-dimensional tensile problem for variational phase-field and gradient-damage models are summarized together with the identification procedure for the phase-field model representative of a prescribed traction-separation cohesive law. The reconstruction of the phase-field model associated with a given cohesive law has been deeply investigated in \cite{Alessi2025b}.
The mechanical response associated with different cohesive laws is thoroughly discussed in \Sec{sec_exp}, highlighting original aspects of the examined models. 
Conclusions and future perspectives are drawn in \Sec{sec_conc}.

\clearpage
\section{The problem setup}
\label{sec_1Dprob}

\begin{figure}[h!]
 \centering
  \begin{subfigure}{\linewidth}
    \fontsize{8pt}{9.6pt}\selectfont 
    \centering
      \includeinkscape[scale=0.8]{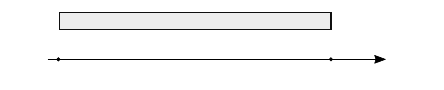}
    \caption{}
    \label{fig_bar} 
  \end{subfigure}
  \\[10mm]
  \begin{subfigure}[b]{0.47\linewidth}
    \centering 
    \fontsize{8pt}{9.6pt}\selectfont 
    \includeinkscape[scale=0.8]{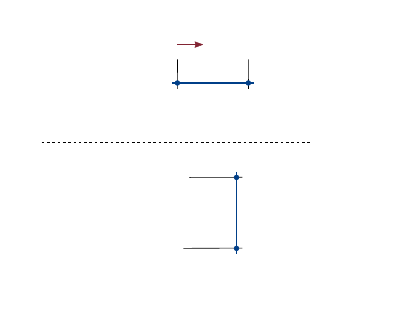}
    \caption{}
    \label{figsub_F} 
  \end{subfigure}
  \begin{subfigure}[t]{0.47\linewidth}
    \centering
    \fontsize{8pt}{9.6pt}\selectfont 
    \vspace{-52mm}
    \includeinkscape[scale=0.8]{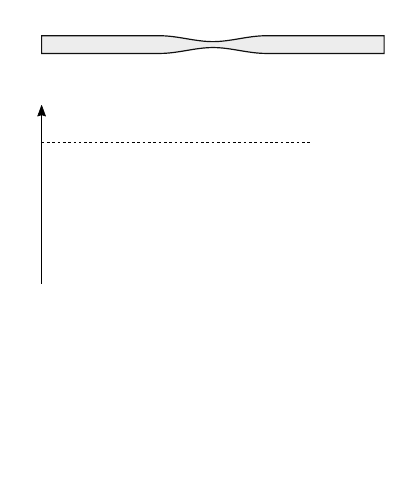}
    \caption{}
    \label{figsub_P} 
  \end{subfigure}
  \caption{The one-dimensional tensile problem (\subref{fig_bar}):  An example of evolution with a sharp-crack (\subref{figsub_F}) and the corresponding phase-field regularization (\subref{figsub_P}).}
  \label{fig_1D}
\end{figure}

Throughout this work, we consider the standard one-dimensional traction bar problem of \Fig{fig_bar}. The bar, initially sound, is assumed to have reference length $L$ and to be made of a homogeneous material with a linear behaviour in the elastic regime and a stress-softening response in the inelastic regime (brittle material). 
The left end ($x = 0$) is kept fixed, while on the right end ($x = L$) a monotonically increasing displacement $U$ is imposed by a hard device, which induces tensile states. Consequently, only positive crack-opening displacements are considered ($\sope \geq 0$).

This paradigmatic problem allows us to illustrate, with closed-form results, the main mechanical features of the phase-field cohesive fracture response of models associated with different traction-separation laws, as discussed at the end of \citep[Sec.~3]{Alessi2025b}. For each considered model we highlight the evolution of key mechanical quantities, including phase-field and displacement profiles as shown in \Fig{figsub_P}.

The extension of the model to higher-dimensional problems is certainly of interest and feasible, although it requires to deal with technical aspects of the formulation beyond the scope of this work and marginally concerning with the identification of the phase-field material functions describing cohesive crack. Indeed, significant contributions in this regard have already been provided in \cite{Focardi2001,Conti2024,Colasanto2025,Conti2025}.

\subsection{The mathematical and engineering frameworks}

In this section we would like to establish a link between the mathematical setting of the first two papers of this work, \cite{Alessi2025a,Alessi2025b} and the present more engineering-oriented contribution.

\paragraph{Mathematical framework}
From a mathematical perspective, as adopted in the main body of this work (Parts~I and~II), the general energy functional \citep[equ. $($1.7$)$]{Alessi2025a} specified for the considered one-dimensional problem of \Fig{fig_bar} reads
\begin{equation}
  \Functeps(u,v,[0,1]):= \int_0^1 \left( \varphi(\eps f^2(v)) |u'|^2 + \frac{\PotDann(1-v)}{4\eps} + \eps |v'|^2 \right) \dx\,
 \label{def_enfuncM}
\end{equation}
with $u(0) = 0$ and $u(1)= U/L$.
In~\eqref{def_enfuncM}, $u$ and $v$ are respectively the displacement and phase-field variable fields on the normalized support~$[0,1]$, with $v$ decreasing from $1$~to~$0$ as the crack evolves, $\varphi(\eps f^2(v))$ a generic degradation function, encompassing the elastic properties of the material, with $\varphi$~satisfying hypotheses (Hp~3-4) and $f(v)\coloneqq\left( \hatf(v)/Q(1-v)\right)^{\sfrac12}\,$ satisfying (Hp~1), $\PotDann(1-v)$ the local phase-field dissipation potential satisfying hypothesis (Hp~1) and $\eps$ an internal length. The hypotheses (Hp~1-4) are reported in \citep[Sec.~1.2]{Alessi2025a}.
The energy functional~\eqref{def_enfuncM} includes, as particular cases, the elastic degradation functions $\varphi(\eps f^2(v))=\widetilde{f}_\eps^2(v)$ and $\varphi(\eps f^2(v))={f}_\eps^2(v)$, with 
\begin{equation}\label{e:f eps tilde}
  \widetilde{f}_\eps(v):=1\wedge \eps^{\sfrac12}f(v)
\end{equation}
and 
\begin{equation}\label{e:f eps}
  f_\eps(v):=\left( \frac{\eps \hatf(v)}{\eps \hatf(v) + Q(1-v)}\right)^{\sfrac12}\,.
\end{equation}
For instance, the degradation function \eqref{e:f eps tilde} has been considered by~\cite{Freddi2017,Lammen2023,Lammen2025}, whereas \eqref{e:f eps} has been considered by~\cite{Wu2017,Wu2018b,Feng2021} with specific choices for the functions $\hatf$ and $Q$. 
The main $\Gamma$-convergence result, which rigorously assess the convergence of the energy functional and minimizers of the phase-field functional \eqref{def_enfuncM} towards the energy functional and minimizers of the free-discontinuity cohesive fracture problem, is provided in the first part of this work, specifically in \citep[Theorem~1.1]{Alessi2025a}.

Throughout this Part~III paper, and according to the motivations illustrated in \Sec{sec_bac_and_mot}, we restrict our attention to degradation functions of the form~\eqref{e:f eps}.
For further details on the phase-field model within the mathematical framework, the reader is invited to refer also to \citep[Sec.~1.2]{Alessi2025a}.

\paragraph{Engineering framework}
From an engineering and mechanical perspective, the emphasis is on the physical interpretation of variables, functions, and parameters, with a focus on material behavior and practical applications. This viewpoint, widely adopted in the engineering community following \cite{Pham2010a,Pham2010c,Marigo2016}, typically expresses the cohesive phase-field energy functional for the problem in \Fig{figsub_P} as
\begin{equation}
  \Enint(\sdisp,\sdam) = \int_0^L \biggl( \,\underbrace{\coh \fDeg(\sdam) E\, \seps^2 \vphantom{\bra{\dfrac{\Gc}{\cw}}} }_{\text{\small $\eqqcolon\enfree(\seps,\sdam)$} }
  + \underbrace{\Gc\bra{ \dfrac{\fDam(\sdam)}{\sintL} + \sintL\, (\sdam')^2 }}_{\text{\small $\eqqcolon\Gc\enfrac(\sdam,\sdam')$} } \, \biggr) \diff x\,,
  \label{def_enfuncE} 
\end{equation} 
where $\sdisp$, $\seps\coloneqq\sdisp'$, and $\sdam$ are, respectively, the displacement, the total strain, and the phase-field variable fields with $\sdam$ increasing from $0$~to~$1$ as the crack evolves, $\fDeg(\sdam)$ is the decreasing material degradation function with the subscript emphasizing its key dependence on the material internal length~$\sintL$ that enables cohesive fracture responses, $E$ denotes the one-dimensional axial stiffness (Young's modulus times the cross section area), $\Gc$ is the one-dimensional fracture toughness and $\fDam(\sdam)$ is the local phase-field dissipation function, an increasing function with respect to~$\sdam$. 
The first term in the integrand, $\enfree(\seps,\sdam)$, represents the internal potential energy whereas the second term represents the dissipated surface fracture energy density with $\enfrac(\sdam,\sdam')$ being the crack density function \cite{Feng2021}.
Their sum, integrated over the entire domain, then represents the total internal energy.
Mechanically, the phase-field variable can be considered as a damage variable and the phase-field model as a gradient-damage model,~\cite{Marigo2016}.
In this view, the following irreversibility condition is often assumed:
\begin{equation}
  \sdamdt \geq 0\,,
\label{def_irr} 
\end{equation} 
aiming at ensuring that the dissipated fracture energy is never decreasing in order to fulfill the second law of thermodynamics and preventing healing effects.
More details on the irreversibility condition are discussed in~\Sec{sec_1Dres}.
For a fully developed crack, the integral over the domain of the crack density function must be equal to one. This requirement is reflected into the phase-field function $\fDam(\sdam)$ by the following condition
\begin{equation}
  4 \int_0^1 \sqrt{\fDam(\sdam)} \diff\sdam = 1\,.
  \label{con_crack} 
\end{equation} 
In standard gradient-damage models, where an arbitrary damage function is usually assumed, condition \eqref{con_crack} is ensured by dividing the fracture toughness by the constant $\cw \coloneqq 4 \int_0^1 \sqrt{\fDam(\sdam)}$.

\paragraph{Link between the mathematical and engineering frameworks}

The displacement variable $u$ has the same meaning in both perspectives, apart from its normalization with respect to the bar length $L$, whereas for the phase-field variable, the following relation holds:
\begin{equation}
 \sdam = 1-v.
 \label{eq_pfcv} 
\end{equation}
By normalizing lengths by $L$ and dividing by $E \,L/2$ in~\eqref{def_enfuncE}, a link between the energy functionals \eqref{def_enfuncM} and \eqref{def_enfuncE} is established. The local phase-field potentials are related by
\begin{equation} \label{eq1}
  \PotDann(1-v) = 4\, k^2 \, \fDam(\sdam)
\end{equation} 
and the internal lengths are related by
\begin{equation}
  \eps = k \dfrac{\ell}{L}\,,
\end{equation}
with $k = 2\Gc\,L / E$ being a material constant.
The degradation function in \eqref{def_enfuncE} consequently reads
\begin{equation}  \label{def_degE} 
  \fDeg(\sdam) = \dfrac{1}{1  + \frac{E}{2 \Gc \ell}\frac{Q(\sdam)}{\lDam(\sdam)}} \, ,
\end{equation} 
with 
\begin{equation} \label{eq2}
  \hatf(v) = \lDam(\sdam).
\end{equation} 
Condition \eqref{con_crack} in terms of $\PotDann$ then becomes
\begin{equation}
  \int_0^1 \sqrt{\PotDann(v)} \diff v= \dfrac{k}{2} \, .
  \label{con_crack2} 
\end{equation}

\subsection{Features of the mechanical model: 1D analytical results}
\label{sec_1Dres} 

In this section, we summarize well-established results and key quantities concerning the mechanical response of the one-dimensional tensile problem for variational phase-field and gradient-damage models, as presented in \cite{Pham2011c, Marigo2016, Wu2017, Chen2021a}, to which we refer the reader for further details.
These results arise from a variational analysis within the energetic formulation framework, where the evolution of a rate-independent system is assumed to be governed by two fundamental principles: a stability condition and an energy balance~\cite{Mielke2006,Mielke2015}.
Most of these results are parameterized with respect to the maximum phase-field state $\sdammax\in\sbra{0,1}$ within the bar.
Throughout this work, we assume without loss of generality that a single localization occurs at $\bar{x} = L/2$, with $\sdam(\bar{x})=\sdammax$. The localization is assumed to remain entirely within the bar during its evolution, avoiding interactions with boundaries. 
This is commonly achieved by enforcing the phase-field variable to vanish at the boundaries. Mechanically, such localization behavior is typically ensured by the design of standard uniaxial test specimens (e.g., dog-bone shapes), in which the ends are thickened and widened to prevent cracking in the grip regions.

\paragraph{Equilibrium, constitutive equation, and phase-field evolution laws}
A first consequence of the variational energetic analysis is that the uniaxial stress is always constant along the bar due to equilibrium, that is,
\begin{equation}
  \sstr' = 0 \qquad \forall \spos \in \bra{0,L}\,,
  \label{def_equ} 
\end{equation}  
where $\sstr$ is linked to the total strain by the constitutive equation
\begin{equation}
  \sstr = \fDeg(\sdam) E \seps \,. 
  \label{def_se} 
\end{equation} 
Another consequence is the set of necessary Kuhn-Tucker conditions in the bulk for the phase-field evolution:
\begin{subequations}
\begin{align}
  \coh \fDeg'(\sdam) E\, \seps^2 + \Gc\bra{ \dfrac{\fDam'(\sdam)}{\sintL} - 2\sintL\, \sdam'' } \geq 0 \,,
  \label{def_fos_con} 
  \\
  \sdamdt \geq 0\,,
  \\
  \bra{\coh \fDeg'(\sdam) E\, \seps^2 + \Gc\bra{ \dfrac{\fDam'(\sdam)}{\sintL} - 2\sintL\, \sdam'' }} \sdamdt = 0\,,
  \label{def_foe_con} 
\end{align}
\label{def_KKT} 
\end{subequations}
$\forall \spos \in \bra{0,L}$.
Equations~\eqref{def_equ} and \eqref{def_KKT} serve as the fundamental basis for deriving the results hereafter presented and used to construct the mechanical responses in \Sec{sec_exp}.

\paragraph{Limit and critical stresses}

The \emph{limit stress} is the stress level associated with the evolving phase-field variable. 
Assuming a stress-softening response, the limit stress is a non-increasing function of $\sdammax$, and reads
\begin{equation}
  \sstr(\sdammax) \coloneqq \sqrt{K\dfrac{\fDam(\sdammax)}{\fCom(\sdammax)}}\,,
  \label{def_strdam} 
\end{equation} 
with $K = 2E\,\Gc/\ell$ being a material constant and $\fCom(\sdam) = 1/\fDeg(\sdam)-1$ being the complementary degradation function. 

Throughout this work, the existence of the following limit is assumed, as stated in \cite[Hp 2]{Alessi2025a}:
\begin{equation}
  \scristr \coloneqq \lim_{\sdam\to 0^{-}} \sstr(\sdammax),
  \label{def_critstr} 
\end{equation} 
with $\scristr$ being the \emph{critical stress}, that is the stress at which the phase-field evolution is triggered and marking the end of the initial elastic response.
The corresponding critical strain then reads $\scrieps \coloneqq \scristr/E$. 
The limit value $\FailureS$ of the mathematical model introduced in \citep[eq. $($1.2$)$]{Alessi2025b} is linked to the critical stress by the following relation
\begin{equation}
  \FailureS = \dfrac{2L}{E} \scristr \, .
  \label{rel_stress} 
\end{equation}

\paragraph{Phase-field profile evolution} 

The symmetric phase-field profiles are described by the following implicit relation:
\begin{equation}
  \abs{x - \bar{x}} = \sintL \int_{\sdam}^{\sdammax} \sqrt{\frac{\fCom(\sdammax)}{\fCom(\sdammax)\fDam(\beta) - \fDam(\sdammax)\fCom(\beta)}}\diff \beta
   \label{eq_damprof} 
\end{equation} 
which scales linearly with the internal length $\sintL$.
Accordingly, half-width of the support of the phase-field localization is given by 
\begin{equation} \label{eq_D} 
  D(\sdammax) = \sintL \int_{0}^{\sdammax} \sqrt{\frac{\fCom(\sdammax)}{\fCom(\sdammax)\fDam(\beta) - \fDam(\sdammax)\fCom(\beta)}}\diff \beta \, .
\end{equation}

\paragraph{Displacement jump (crack opening)} 

For a one-dimensional tensile problem, the displacement jump amplitude associated with the evolving cohesive crack, that is the crack opening, can easily be determined.
From the constitutive relation~\eqref{def_se}, the total strain reads
\begin{equation}
  \seps = \dfrac{\sstr}{\fDeg(\sdam)\,E} = \dfrac{\sstr}{E}\bra{1 + \fCom(\sdam)}\, .
  \label{def_strain}
\end{equation} 
Integrating \eqref{def_strain} over the entire bar, the right-end global displacement is then given by 
\begin{equation}
  U = \dfrac{\sstr}{E}\int_{0}^{L} \bra{1 + \fCom(\sdam(\spos))} \diff \spos
    = \dfrac{\sstr}{E}L + \underbrace{2\dfrac{\sstr}{E}\int_{\bar\spos}^{\bar\spos+D} \fCom(\sdam(\spos)) \diff \spos}_{\text{\small $\eqqcolon \delta$ \newline (crack opening)} }  \,.
\label{eq_U} 
\end{equation} 
The first term in the last expression clearly represents the elastic contribution to the overall displacement. Consequently, the second term corresponds to the crack displacement jump $\sope$. By changing the variable $x\to\sdam$ and using \eqref{def_critstr}, the displacement jump as a function of the maximum phase-field value is given by
\begin{equation}
\begin{aligned}
  \sope(\sdammax) 
    &= \dfrac{2\,\sstr}{E}\int_{0}^{D} \fCom(\sdam(\spos-\bar\spos)) \diff \spos \\
    &= \dfrac{2\,\sstr \, \sintL}{E} \int_{0}^{\sdammax} \sqrt{\frac{\fCom(\sdammax)}{\fCom(\sdammax)\fDam(\beta) - \fDam(\sdammax)\fCom(\beta)}} \fCom(\beta) \diff \beta \\
    &= \dfrac{4\,\Gc}{\cw \scristr} \sqrt{\dfrac{\fDam'(0)}{\fCom'(0)}} 
    \int_{0}^{\sdammax} \sqrt{\frac{\fDam(\sdammax)}{\fCom(\sdammax)\fDam(\beta) - \fDam(\sdammax)\fCom(\beta)}} \fCom(\beta) \diff \beta  \,.
\label{def_sope} 
\end{aligned}
\end{equation} 
The ultimate crack opening, that may be finite or infinite, is given by
\begin{equation}
  \scriope \coloneqq \lim_{\sdam \to 1} \sope(\sdam)  \,. 
  \label{def_scriope} 
\end{equation}

\paragraph{Traction-separation cohesive law}
Combining \eqref{def_strdam} and \eqref{def_sope} yields the traction-separation law parameterized with respect to the maximum phase-field state:
\begin{equation} \label{def_param} 
  \sstr(\sope) \,:\, \sdammax \mapsto \bra{\sope(\sdammax),\sstr(\sdammax)}  \,.
\end{equation}
Notably, neither $\sstr$ nor $\sope$ depend on $\sintL$, unlike models based on~\eqref{e:f eps tilde}.

In \citep[Sec. 3.3]{Alessi2025b}, where the construction of a phase-field model associated with a specific cohesive law is explored, the assumption $Q = \PotDann/\FailureS$ is made throughout, without loss of generality. By applying this assumption to the mechanical model using the relations \eqref{eq1} and \eqref{rel_stress}, the degradation function \eqref{def_degE} becomes:
\begin{equation} \label{def_degE2} 
  \fDeg(\sdam) = \dfrac{1}{1  + \frac{2\,\Gc\,E}{\ell\,\scristr^2}\frac{\fDam(\sdam)}{\lDam(\sdam)}} \,.
\end{equation} 
Then, \eqref{def_degE2} explicitly shows how \eqref{def_enfuncE} depends on the critical stress $\scristr$ and the fracture toughness $\Gc$, two material parameters with clear mechanical significance that characterize the failure behavior and can be easily identified through standard experiments.

\paragraph{Dissipated fracture energy}

The fracture toughness $\Gc$ is equal to the area of the region of the $(\sope,\sstr)$ plane delimited by the $\sstr$ axis and the curve $\sstr \mapsto \sope(\sstr)$.
Indeed
\begin{equation}
\begin{aligned}
  \int_{0}^{\scristr} \sope(\sstr) \diff \sstr 
    &= \int_{0}^{\scristr} \dfrac{2\,\sstr \, \sintL}{E} \int_{0}^{\sdammax(\sstr)} \sqrt{\frac{\fCom(\sdammax)}{\fCom(\sdammax)\fDam(\beta) - \fDam(\sdammax)\fCom(\beta)}} \fCom(\beta) \diff \beta \diff \sstr \\
    &=  \dfrac{2\, \sintL}{E} \int_{0}^{1} \int_{0}^{\sstr(\beta)} \frac{\sstr \, \fCom(\beta)}{\sqrt{\fDam(\beta) - \frac{\sstr^2}{K} \fCom(\beta)}}  \diff \sstr  \diff \beta \\
    &= - \dfrac{2\, \sintL\,K}{E} \int_{0}^{1} \sbra{\sqrt{\fDam(\beta) - \frac{\sstr^2}{K} \fCom(\beta)}\:}_{0}^{\sstr(\beta)}  \diff \beta \\
    &= \dfrac{2\, \sintL\,K}{E} \int_{0}^{1} \sqrt{\fDam (\beta)} \diff \beta \\
    &= \Gc \,.
\end{aligned}
\label{eq_Gc} 
\end{equation} 
For all models except Dugdale's, presented in the next section, the fracture toughness corresponds to the dissipated energy.
Indeed, in Dugdale's cohesive phase-field model discussed in \Sec{sec_D}, a snap-back response occurs during the evolution, causing the dissipated energy to be greater than $\Gc$, the dissipated fracture energy. It is worth remarking that even in this case, \eqref{eq_Gc} remains valid.
The surface fracture energy as a function of the crack opening, that is the fracture energy dissipated up to the a given crack opening, will be denoted as $G=G(\sope)$.
Within the mathematical framework, the crack opening, the surface fracture energy and the limit stress are respectively denoted as $s$, $g_0(s)$ and $g'_0(s)$.

\paragraph{Irreversibility}

The issue of modeling irreversibility in phase-field models approximating sharp fracture problems is still an open problem.
When the phase-field variable is treated as a damage variable, the irreversibility condition \eqref{def_irr} is typically imposed \cite{Marigo2016}.
However, phase-field profiles constructed from \eqref{eq_damprof} do not always guarantee the fulfillment of the irreversibility condition \eqref{def_irr}, as discussed in \cite{Wu2024a}.

In sharp fracture models, irreversibility applies to the fracture set $\Gamma$, approximated by the region where $\sdam=1$, ensuring that $\Gamma(t_{i})\subseteq\Gamma(t_{j}), \: \forall t_{i} \leq t_{j}$, where~$t$ is the time-evolution parameter,~\cite{Bourdin2008}.
From this perspective, \eqref{def_irr} may be overly restrictive. Moreover, in cohesive fracture models, irreversibility also influences the unloading behavior in traction-separation laws, an aspect considered even more rarely,~\cite{Crismale2016a,Bonacini2021}.

In this work, we interpret the phase-field variable primarily as a mathematical regularization tool rather than a strict damage variable, and thus do not enforce \eqref{def_irr} on its evolution.
However, this issue will be further discussed in the next section, where different models yielding the same cohesive response suggests the possibility of finding models satisfying \eqref{def_irr}.

\subsection{Identifying the phase-field model associated with a prescribed cohesive fracture law}
\label{sec_ident} 
In Part~II of this work~\cite{Alessi2025b}, we have shown, within a mathematical framework, how to assign a given cohesive law, either linear or superlinear for small amplitude of the displacement jump, by appropriately selecting the parameters of the phase-field model in \eqref{def_enfuncM}.
More precisely, by either fixing the local phase-field dissipation potential $\PotDann$ and choosing an appropriate degradation function~$f_\eps(v)$ or viceversa, we define the corresponding functionals $\Functeps$ and demonstrate that $g_0$ equals the surface energy density of their $\Gamma$-limit~$g$, as established in~\citep[Theorem~1.1]{Alessi2025a}.
The procedure to construct a phase-field model associated with a prescribed cohesive law is developed in detail in Part~II~\cite{Alessi2025b}, under the following assumptions for the linear case
\begin{itemize}
  \item[(Hp~$6$)] $g_0\in C^1([0,\infty))$, $g_0^{-1}(0)=\{0\}$, $g_0$ bounded and non-decreasing;
  \item[(Hp~$7$)] $g_0$ concave,  $g_0'(0)=\FailureS\in(0,\infty)$, $g'_0$ strictly decreasing on $[0,(\sfratt)_0)$ where 
    \[
      (\sfratt)_0:=\sup \{s \in [0,\infty):\, g_0(s)<\lim_{s\to\infty}g_0(s)\}\in(0,\infty]\,.
    \]
  \item[(Hp~$8$)] $R$, the auxiliary function based on $g_0$ defined in \citep[$($3.1$)$]{Alessi2025b}, is convex on $[0,\FailureS^2]$ and $W^{1,p}((0,\FailureS^2))$ for some $p>2$.
\end{itemize}
together with the standard assumptions (Hp~1-4) for the validity of Theorem 1.1.
For the superlinear case, which it not considered in the following examples, assumptions (Hp 6) and (Hp 7) are replaced by similar hypotheses (Hp 6') and (Hp 7'),~\cite{Alessi2025b}.

Operatively, the steps required to deduce the phase-field model corresponding to a given cohesive law are:
\begin{enumerate}
 \item Consider the physical traction-separation law $\sstr(\sope)$ for a given cohesive law;
 \item Define the non-dimensional cohesive law $g_0(s)$, adopted within the mathematical framework in \cite{Alessi2025b}, from the relation
    \begin{equation}
      g'_0(s) \coloneqq \dfrac{\sstr(\kscal\,s)}{\sstr(0)} \,,
      \label{def_gprime} 
    \end{equation} 
    with $g_0(0)=0$ and $\kscal>0$ being a scaling factor for the non-dimensional displacement jump~$s$ which encompasses a characteristic crack opening reference length. Of course, an admissible cohesive law $g_0(s)$ and the corresponding traction-separation law~$g'_0(s)$ must satisfy (Hp~6) and (Hp~7); 
  \item Assuming, without loss of generality, $Q(t)=\omega(t)$, determine the functions $\hatf$ and $\PotDann$, fully characterizing the phase-field regularizing energy functional~\eqref{def_enfuncM}, by adopting the construction in Theorem.~3.1 or 3.2 in \cite{Alessi2025b};
  \item Fix $\kscal$ through \eqref{con_crack2};
  \item Reconstruct the material functions $\fDam$, $\lDam$ and $\fDeg$ of the engineering phase-field energy functional through \eqref{eq1}, \eqref{eq2} and \eqref{def_degE2}, respectively.
\end{enumerate}

It is worth emphasizing that, for a given admissible cohesive law, multiple combinations of $\cbra{\hatf,\PotDann}$ may be deduced. This aspect and its implications will be highlighted throughout the examples in \Sec{sec_exp}.

\section{Examples of the mechanical response of different cohesive laws}
\label{sec_exp} 

In this section, we present the mechanical response of the uniaxial tensile bar test of \Fig{fig_bar} for different cohesive-laws considered in the examples of Section~3 in \cite{Alessi2025b}. The phase-field model associated with each assigned cohesive law is obtained through the construction summarized in \Sec{sec_ident}. Specifically we consider throughout this section the following softening cohesive laws: Linear (\Fig{fig_E_L} and \Sec{sec_L}), bi-linear (\Fig{fig_E_L2} and \Sec{sec_L2}), exponential (Barenblatt like, \Fig{fig_E_E} and \Sec{sec_E}), hyperbolic (Fig.s~\ref{fig_E_H} and~\ref{fig_E_H2}, and \Sec{sec_H}) and Dugdale (\Fig{fig_E_D} and \Sec{sec_D}).

\newcommand{\myfig}[1]{\includegraphics[page=#1,scale=0.8,trim=3mm 0 0 1mm, clip]{ACF_3-Examples}}
\newcommand{\mysubref}[1]{(\subref{#1})}

\begin{figure}[h!]
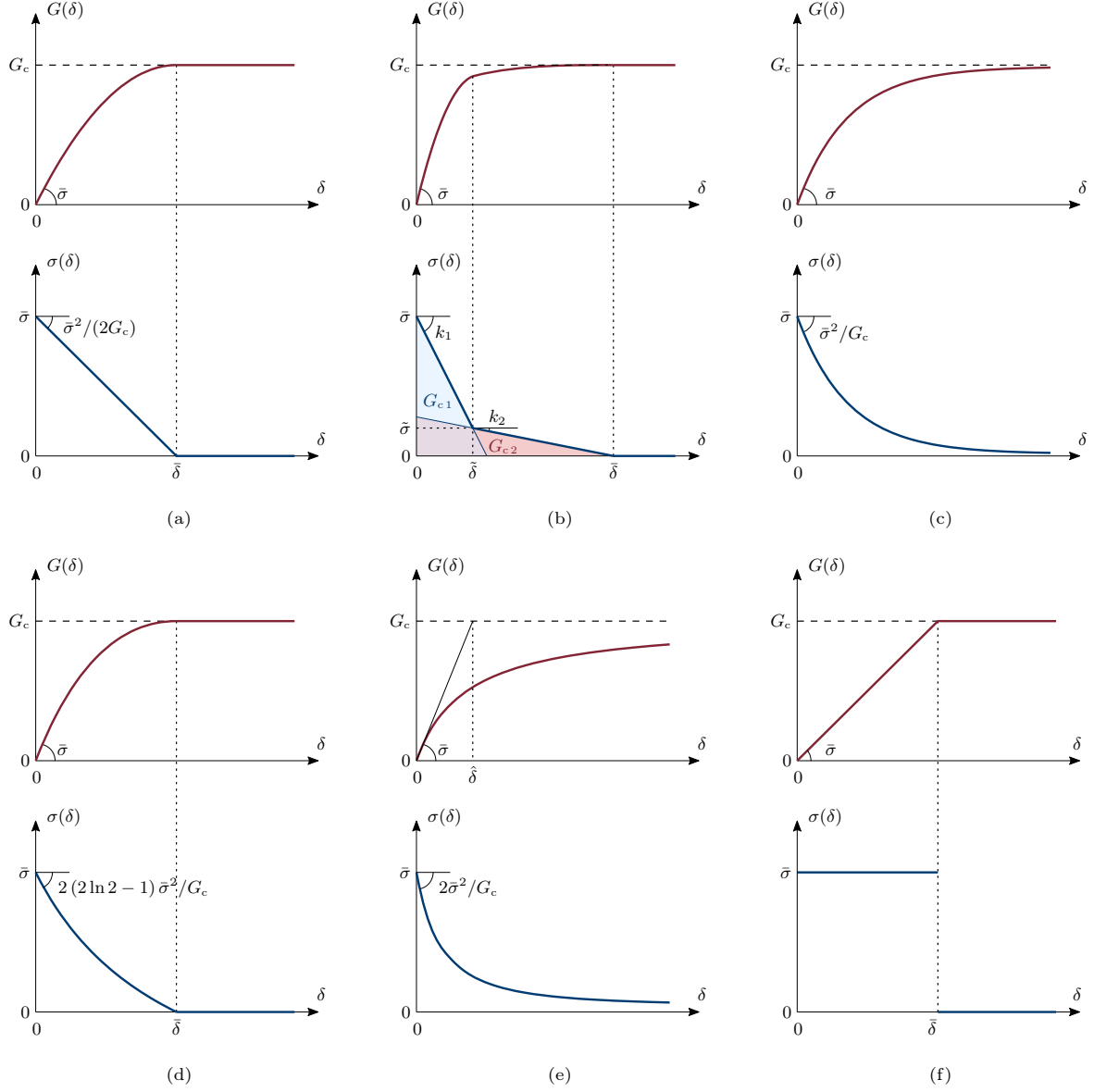

  \centering
  \begin{subfigure}{0.32\linewidth}
    \centering
      \myfig{2}
    \caption{} \label{fig_E_L} 
  \end{subfigure}
  \hfill
  \begin{subfigure}{0.32\linewidth}
    \centering
      \myfig{3}
    \caption{} \label{fig_E_L2} 
  \end{subfigure}
  \hfill
  \begin{subfigure}{0.32\linewidth}
    \centering
      \myfig{6}
    \caption{} \label{fig_E_E} 
  \end{subfigure}
  \\[1ex]
  \begin{subfigure}{0.32\linewidth}
    \centering
      \myfig{4}
    \caption{} \label{fig_E_H} 
  \end{subfigure}
  \hfill
  \begin{subfigure}{0.32\linewidth}
    \centering
      \myfig{5}
    \caption{} \label{fig_E_H2} 
  \end{subfigure}  
  \hfill
  \begin{subfigure}{0.32\linewidth}
    \centering
      \myfig{1}
    \caption{} \label{fig_E_D} 
  \end{subfigure}
\caption{Qualitative trends of the considered softening cohesive laws: \mysubref{fig_E_L} linear, \mysubref{fig_E_L2} bilinear, \mysubref{fig_E_E} exponential (Barenblatt like), hyperbolic (\mysubref{fig_E_H} linear and \mysubref{fig_E_H2} quadratic) and \mysubref{fig_E_D} Dugdale.}
\label{fig_Ex} 
\end{figure}

These laws are widely used by the engineering community for describing the fracture process of materials exhibiting a quasi-brittle failure behaviour as, for instance, concrete, fibre-reinforced composites or geological materials \cite{Guinea1994,Roesler2007a,Park2011,Dourado2012}. Such quasi-brittle failures often involves a combination of progressive micro-cracking and crack bridging or aggregate interlocking. 
The geometric and constitutive parameters adopted throughout this section are reported in \Tab{tab_par} and are the same as in \citep[Sec. 7.3]{Feng2021}.

\begin{table}[h!]
  \centering
  \small
  \renewcommand{\arraystretch}{1.4}
  \begin{tabular}{cr@{\,}l}
    \toprule
    $L$ & 200 & \si{\mm} \\
    $E$ & 30000 & \si{\MPa} \\
    $\Gc$ & 0.12 & \si{\newton\over\mm} \\
    $\scristr$ & 3 & \si{\MPa} \\
    $\sintL$ & \multicolumn{2}{c}{variable} \\
    \bottomrule
  \end{tabular}
  \caption{Geometric and constitutive parameters.}
  \label{tab_par} 
\end{table}

\FloatBarrier
\subsection{Linear softening cohesive law}
\label{sec_L} 

The linear softening cohesive law is completely defined once two of the following three parameters are assigned: the tensile strength $\scristr$, the fracture toughness $\Gc$ and the ultimate crack opening~$\scriope$, \Fig{fig_E_L}. The linear traction-separation law reads
\begin{equation}
  \sstr(\sope) =
  \begin{cases}
    \scristr \bra{1 - \dfrac{\sope}{\scriope}}, \quad &\text{for $\sope \leq \scriope$} \\
    0, \quad &\text{for $\sope > \scriope$}
  \end{cases}  \,,
  \label{eq_TSL_L} 
\end{equation} 
where the ultimate crack opening is given by $\scriope = 2\,\Gc/\scristr$.

Replacing $\delta$ by $k_{\delta}\,s$, with $k_{\delta}$ being a generic constant, and fixing $k_s = k_{\delta}/\scriope$, \eqref{eq_TSL_L} yields, according to \eqref{def_gprime}, the following linear traction-separation law within the mathematical framework 
 \begin{equation}
    g'(s)=
    \begin{cases}
        1- \kscal s, \; & \text{for $s \leq \frac{1}{\kscal}$} \\
        0, \; & \text{for $s > \frac{1}{\kscal}$}
    \end{cases}  \,.
  \label{eq_TSLM_L} 
\end{equation}
Thus, $k_s$ encompasses here the characteristic crack opening length~$\scriope$. A similar rescaling procedure will be adopted to all subsequent examples.

\begin{figure}[h]
\centering
\begin{subfigure}[b]{0.4\linewidth}
  \begin{overpic}[width=0.8\linewidth]{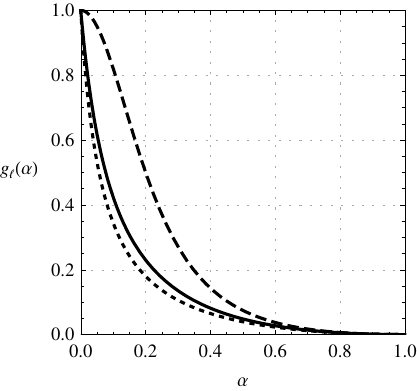}
    \put(59,75){\LLeg}
  \end{overpic}
  \caption{}
  \label{gLPlot} 
\end{subfigure}
\begin{subfigure}[b]{0.4\linewidth}
  \begin{overpic}[width=0.8\linewidth]{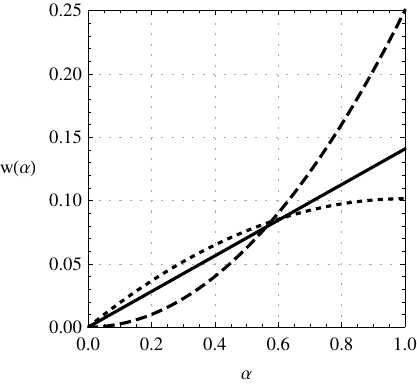}
    \put(22,75){\LLeg}
  \end{overpic}
  \caption{}
  \label{wLPlot}
\end{subfigure}
  \caption{(\subref{gLPlot}) Degradation and (\subref{wLPlot}) local dissipation functions for the linear cohesive fracture models \eqref{mod_L1}, \eqref{mod_L2} and \eqref{mod_L12}.}
  \label{fig_gwL_func} 
\end{figure}

In \cite[Sec. 3.3.2]{Alessi2025b}, three different models, all describing the same traction-separation law~\eqref{eq_TSLM_L}, have been derived, sharing the same crack opening scaling factor $k_s=1/(2k)$. The corresponding engineering material functions $\lDam$ and $\fDam$ for these models read:

\begin{equation}
\begin{cases}
  \lDam^{-1}(t) = 1-\bra{1-\frac{2}{\pi } \bra{\arccos(\sqrt{t}) + \sqrt{t-t^2}}}^{2/3}\\
  \fDam(\sdam) = \dfrac{9}{64}\sdam
\end{cases},
\tag{\texttt{L}$_1$}
\label{mod_L1} 
\end{equation}
\begin{equation}
\begin{cases}
  \lDam^{-1}(t) = 1-\bra{1-\frac{2}{\pi } \bra{\arccos(\sqrt{t}) + \sqrt{t-t^2}}}^{1/2}\\
  \fDam(\sdam) = \dfrac{\sdam^2}{4}
\end{cases},
\tag{\texttt{L}$_2$}
\label{mod_L2} 
\end{equation}
and
\begin{equation}
\begin{cases}
  \lDam(\sdam)= (1-\sdam)^2 \\
  \fDam(\sdam)= \dfrac{2\sdam - \sdam^2}{\pi^2}
\end{cases}.
\tag{\texttt{L}$_{12}$}
\label{mod_L12} 
\end{equation} 
The functions $\lDam(\sdam)$ for models \eqref{mod_L1} and \eqref{mod_L2} do not admit an explicit analytical inverse. As a consequence, interpolation functions to reconstruct $\lDam(\sdam)$ are hereafter adopted.
Models \eqref{mod_L1} and \eqref{mod_L2} have been derived by fixing $\hatf$ and deducing, by means of Theorem~3.1 in~\cite{Alessi2025b}, $\PotDann$ whereas, for model \eqref{mod_L12}, $\PotDann$ has been first fixed and (the inverse of) $\hatf$ deduced by means of~Theorem~3.2 in~\cite{Alessi2025b}. Model \eqref{mod_L12} corresponds to the linear model considered in \citep[Section 4.1]{Feng2021}.
The degradation function $\fDeg$ and the local material dissipation function $\fDam$ trends for the three linear models are depicted in \Fig{fig_gwL_func}.
A major difference between these three models is the function behaviour of $\fDam$, namely: linear for \eqref{mod_L1}, convex for \eqref{mod_L2} and concave for \eqref{mod_L12}. 
It is worth also noting that the local phase-field dissipation potential~$\fDam$ for \eqref{mod_L2}, which is essentially the same of the \texttt{AT}$_2$ phase-field model for brittle fracture~\cite{Tanne2018}, lacks of a linear term. Nevertheless, this model is still able to describe an explicit elastic response associated with a non-vanishing critical stress. Therefore, the presence of a linear term in the local material dissipation function is not a necessary but only sufficient condition for observing an explicit elastic response.

\begin{figure}[h]
\centering
\begin{subfigure}[b]{0.31\linewidth}
  \begin{overpic}[width=\linewidth]{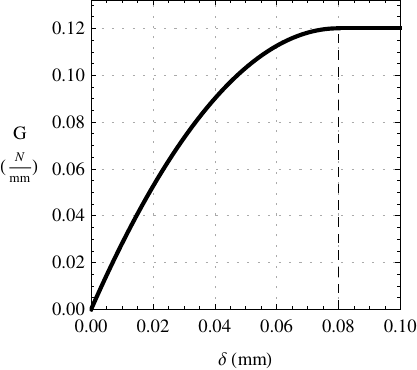}
        \put(60,25){\LLeg}
  \end{overpic}
  \caption{}
  \label{GLPlot} 
\end{subfigure}
\quad
\begin{subfigure}[b]{0.31\linewidth}
  \begin{overpic}[width=\linewidth]{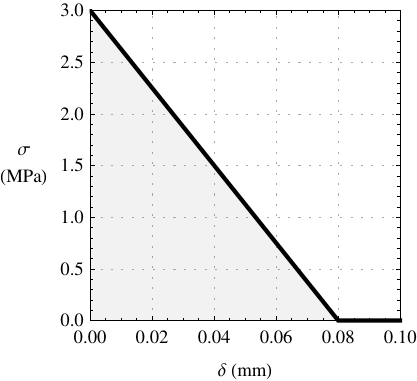}
        \put(35,32){$\Gc$}
        \put(58,75){\LLeg}
  \end{overpic}
  \caption{}
  \label{tslawLPlot1} 
\end{subfigure}
\\[2ex]
\begin{subfigure}[b]{0.31\linewidth}
  \begin{overpic}[width=\linewidth]{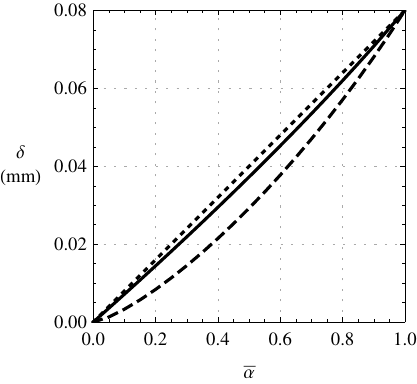}
    \put(21,75){\LLeg}
  \end{overpic}
  \caption{}
  \label{tslawLPlot2} 
\end{subfigure}
\quad
\begin{subfigure}[b]{0.31\linewidth}
  \begin{overpic}[width=\linewidth]{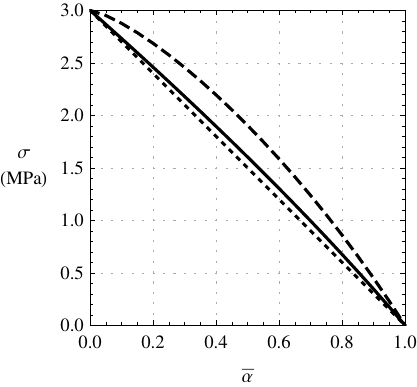}
    \put(58,75){\LLeg}
  \end{overpic}
  \caption{}
  \label{tslawLPlot3} 
\end{subfigure}
  \caption{(\subref{GLPlot}) Surface fracture energy and (\subref{tslawLPlot1}) traction-separation law as a function of the displacement jump, obtained using the parametrisation \eqref{def_param} of the displacement jump (\subref{tslawLPlot2}) and the critical stress (\subref{tslawLPlot3}) with respect to the maximum phase-field value for all considered linear models \eqref{mod_L1}–\eqref{mod_L12}.}
  \label{fig_GL} 
\end{figure}

\begin{figure}[h]
\centering
  \begin{overpic}[scale=0.8]{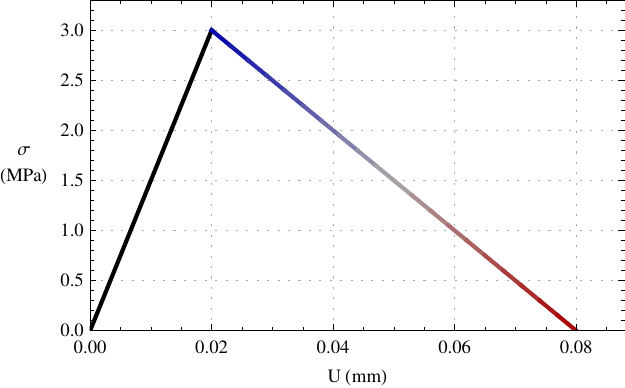}
    \put(76,55){\damBar[$\sdammax$]}
  \end{overpic}
  \caption{Global response of the bar, exhibiting an independent response both with respect to the linear models (\eqref{mod_L1}, \eqref{mod_L2} and \eqref{mod_L12}) and the different internal lengths~\eqref{def_lL}. The initial elastic response is shown in black. The dependence on the maximum phase-field state is also highlighted. }
  \label{fig_UFL} 
\end{figure}

Fig.s \ref{GLPlot} and \ref{tslawLPlot1} respectively represent the surface fracture energy and the traction-separation law as a function of the displacement jump for all the considered linear models. 
As also evident from the global response in \Fig{fig_UFL}, which exhibits a linear elastic response followed by a linear softening damage response, the tree models \eqref{mod_L1}-\eqref{mod_L12} are describing the same cohesive fracture response.

\begin{figure}[h]
\centering
\begin{subfigure}[b]{\linewidth}
  \centering
  \begin{overpic}[width=0.8\linewidth]{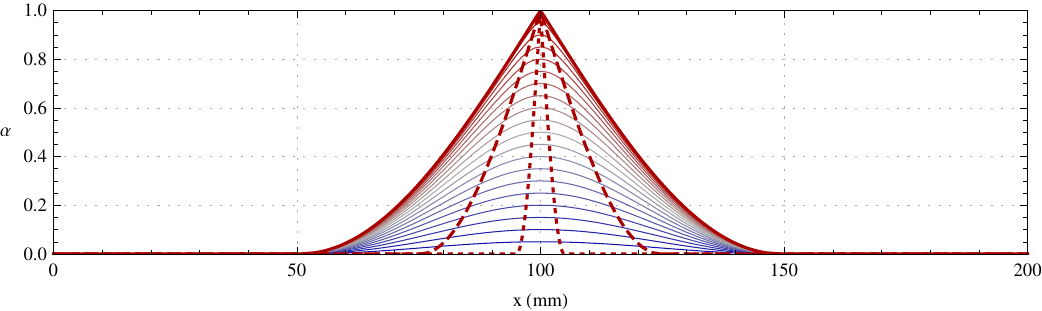}
    \put(81,22){\damBar}
    \put(76,12){\lLeg}
  \end{overpic}
  \caption{}
  \label{fig_DamProf1} 
\end{subfigure}
\\[2ex]
\begin{subfigure}[b]{\linewidth}
  \centering
  \begin{overpic}[width=0.8\linewidth]{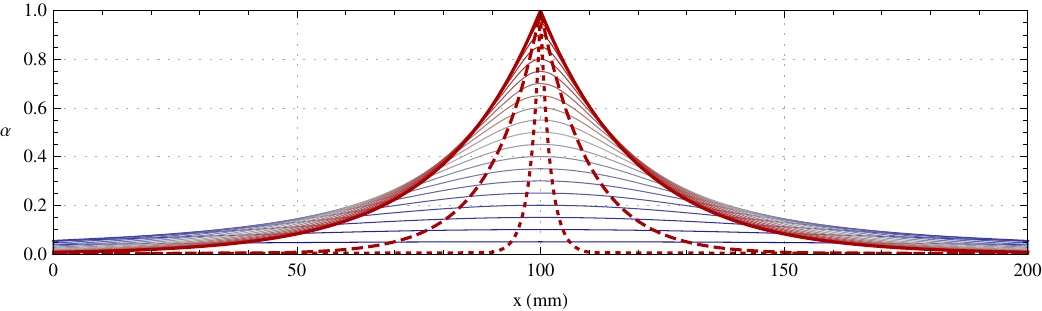}
    \put(82,24){\damBar}
    \put(77,14){\lLeg}
  \end{overpic}
  \caption{}
  \label{fig_DamProf2} 
\end{subfigure}
\\[2ex]
\begin{subfigure}[b]{\linewidth}
  \centering
  \begin{overpic}[width=0.8\linewidth]{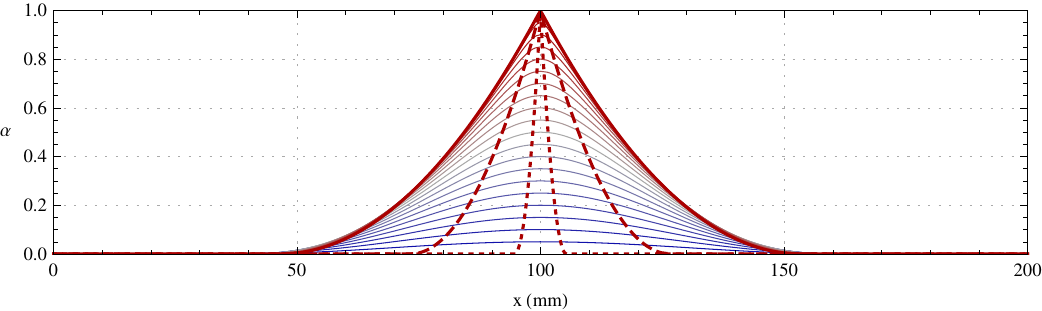}
    \put(81,22){\damBar}
    \put(76,12){\lLeg}
  \end{overpic}
  \caption{}
  \label{fig_DamProf3} 
\end{subfigure}
  \caption{Phase-field profiles at given maximum phase-field states for the linear models: (\subref{fig_DamProf1}) for \eqref{mod_L1}, (\subref{fig_DamProf2}) for \eqref{mod_L2} and (\subref{fig_DamProf3}) for \eqref{mod_L12}. The ultimate phase-field profiles associated with the three internal lengths~\eqref{def_lL} are also highlighted.}
  \label{fig_DamProfL} 
\end{figure}

\begin{figure}[h]
\centering
\footnotesize
  \begin{overpic}[scale=0.9]{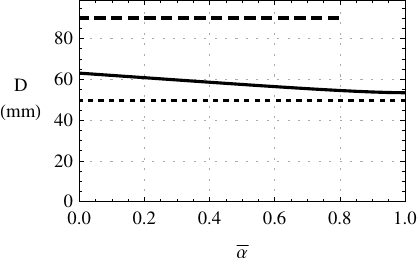}
    \put(83,56){$\uparrow +\infty$}
    \put(68,23){\LLeg}
  \end{overpic}
  \caption{Half-width of the support of the phase-field profiles for the three linear models \eqref{mod_L1}, \eqref{mod_L2} and \eqref{mod_L12} with respect to the fracture evolution ($\sintL=\SI{10}{mm}$).}
  \label{fig_DL} 
\end{figure}

To highlight the role of~$\sintL$, three different values for the internal length have been considered, namely:
\begin{equation} \label{def_lL} 
  \sintL = \cbra{10,5,1} \quad \text{(mm)} \,.
\end{equation}
It turns out, as predicted, that the global response is independent of the value of the internal-length, provided it is sufficiently small compared to the domain size, \Fig{fig_UFL}. The internal length has instead a significant impact on the phase-field localization size, as evident in \Fig{fig_DamProfL}, which reports the phase-field profiles evolutions, obtained from \eqref{eq_damprof}, for the linear models together with the ultimate phase-field profiles associated with the three internal lengths~\eqref{def_lL}.
The trends of the half-width localization support~\eqref{eq_D} with respect to the fracture evolution are highlighted in \Fig{fig_DL}.
It is then immediately clear that although the three linear models describe the same cohesive fracture response and global response, the behaviour of their corresponding phase-field profiles is very different. For instance, the phase-field profile for \eqref{mod_L12} evolves with a constant and finite support, therefore automatically satisfying the irreversibility condition \eqref{def_irr}, whereas for \eqref{mod_L2} the support in unbounded. Clearly, for these examples, the phase-field at the boundaries has been allowed to assume values different from zero. The phase-field profiles for model \eqref{mod_L1} instead have a finite support that slightly shrinks during the evolution.

\begin{figure}[h]
\centering
  \begin{overpic}[width=0.8\linewidth]{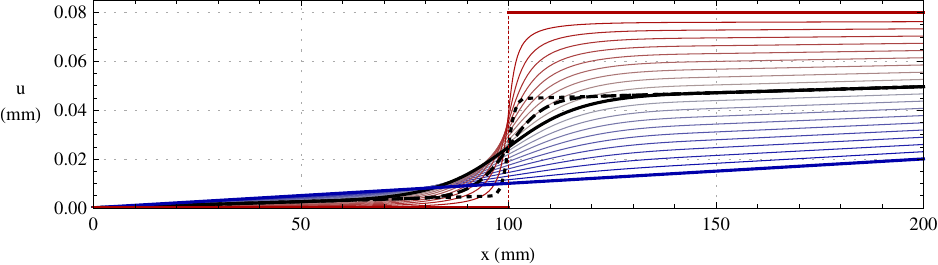}
    \put(10,24){\damBar}
    \put(10,14){\lLeg[0.5]}
  \end{overpic}
  \caption{Displacement profiles at given maximum phase-field states for \eqref{mod_L12} and $\sintL = \SI{10}{mm}$. The displacement profiles at $\sdammax = 0.5$ for smaller internal lengths are also highlighted.}
  \label{fig_UProfL} 
\end{figure}

We also report the evolution of the displacement profile for \eqref{mod_L12} where the convergence of the smooth transition towards a sharp is noticeable as the internal length decreases,~\Fig{fig_UProfL}.

\FloatBarrier

\subsection{Bilinear softening cohesive law}
\label{sec_L2} 

\begin{figure}[h]
\centering
\begin{subfigure}[b]{0.4\linewidth}
  \includegraphics[width=0.8\linewidth]{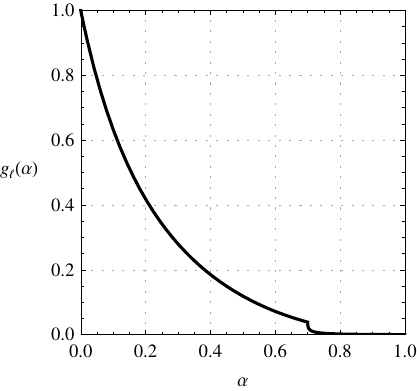}
  \caption{}
  \label{gBPlot} 
\end{subfigure}
\begin{subfigure}[b]{0.4\linewidth}
  \includegraphics[width=0.8\linewidth]{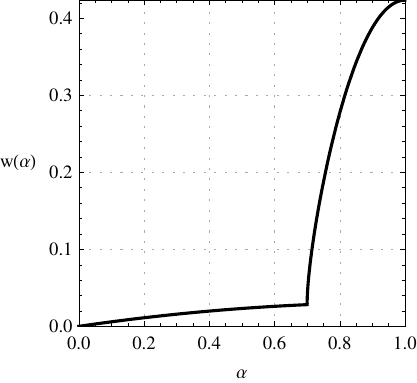}
  \caption{}
  \label{wBPlot} 
\end{subfigure}
  \caption{(\subref{gBPlot}) Degradation and (\subref{wBPlot}) local dissipation functions for the bilinear cohesive fracture model~\eqref{mod_B}.}
  \label{fig_gwB_func} 
\end{figure}

\begin{figure}[h]
\centering 
  \begin{subfigure}[b]{0.31\linewidth}
    \includegraphics[width=\linewidth]{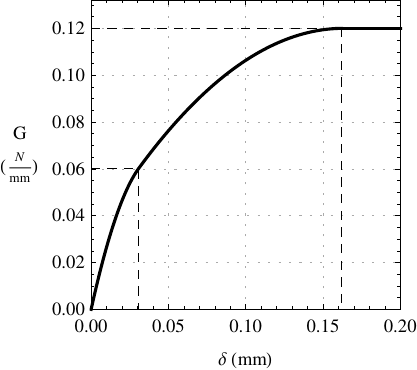}
    \caption{}
    \label{GBPlot} 
  \end{subfigure}
  \quad
  \begin{subfigure}[b]{0.31\linewidth}
    \begin{overpic}[width=\linewidth]{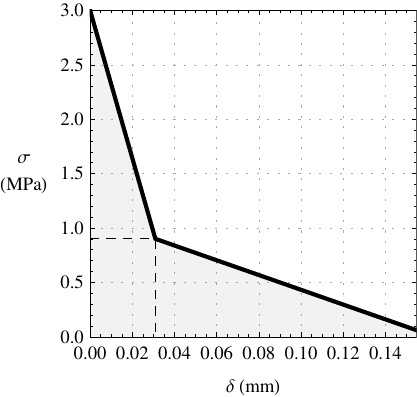}
      \put(26,25){$\Gc$}
    \end{overpic}
    \caption{}
    \label{tslawBPlot1} 
  \end{subfigure}
\\[2ex]
  \begin{subfigure}[b]{0.31\linewidth}
    \includegraphics[width=\linewidth]{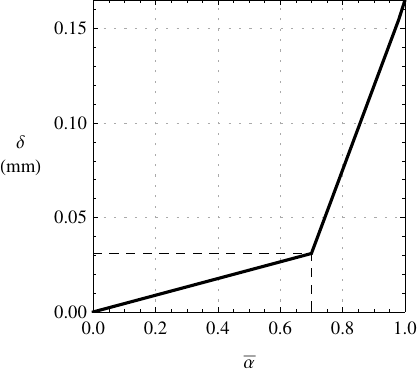}
    \caption{}
    \label{tslawBPlot2} 
  \end{subfigure}
  \quad
  \begin{subfigure}[b]{0.31\linewidth}
    \includegraphics[width=\linewidth]{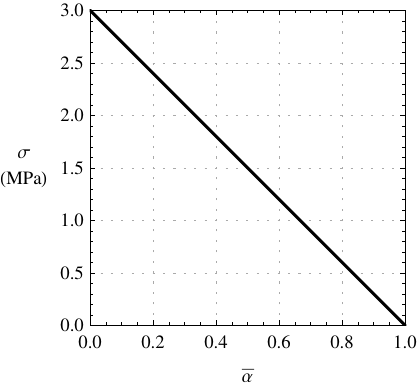}
    \caption{}
    \label{tslawBPlot3} 
  \end{subfigure}
  \caption{(\subref{GBPlot}) Surface fracture energy and (\subref{tslawBPlot1}) traction-separation law of the bilinear model~\eqref{mod_B} as a function of the displacement jump, obtained using the parametrization \eqref{def_param} of the displacement jump (\subref{tslawBPlot2}) and the critical stress~(\subref{tslawBPlot3}) with respect to the maximum phase-field value.}
  \label{fig_GB} 
\end{figure}

\begin{figure}[h]
\centering
  \begin{overpic}[scale=0.8]{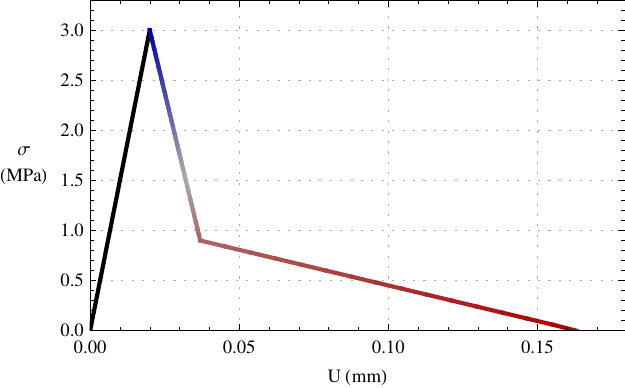}
    \put(76,55){\damBar[$\sdammax$]}
  \end{overpic}
  \caption{Global response for the bilinear model~\eqref{mod_B}, independent of the internal length. The initial elastic response is shown in black. The dependence on the maximum phase-field state is also highlighted.}
  \label{fig_UFB} 
\end{figure}

\begin{figure}[h]
\centering
  \begin{overpic}[width=0.8\linewidth]{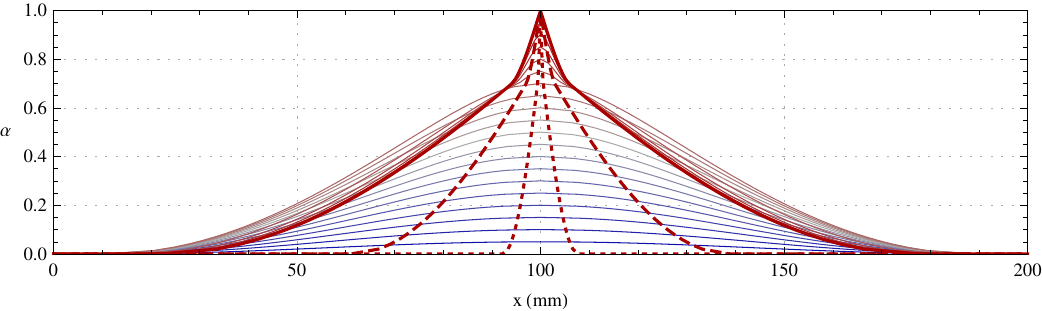}
    \put(81,24){\damBar}
    \put(76,15){\lLeg}
  \end{overpic}
  \caption{Phase-field profiles at given maximum phase-field states for the bilinear model~\eqref{mod_B}. The ultimate phase-field profiles associated with the three internal lengths~\eqref{def_lL} are also highlighted.}
  \label{fig_DamProfB} 
\end{figure}

\begin{figure}[h]
\centering
\footnotesize
  \includegraphics[scale=0.9]{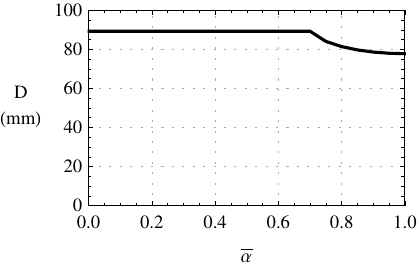}
  \caption{Half-width of the support of the phase-field profile for the bilinear model~\eqref{mod_B} with respect to the fracture evolution ($\sintL=\SI{10}{mm}$).}
  \label{fig_DB} 
\end{figure}

\begin{figure}[h]
\centering
  \begin{overpic}[width=0.8\linewidth]{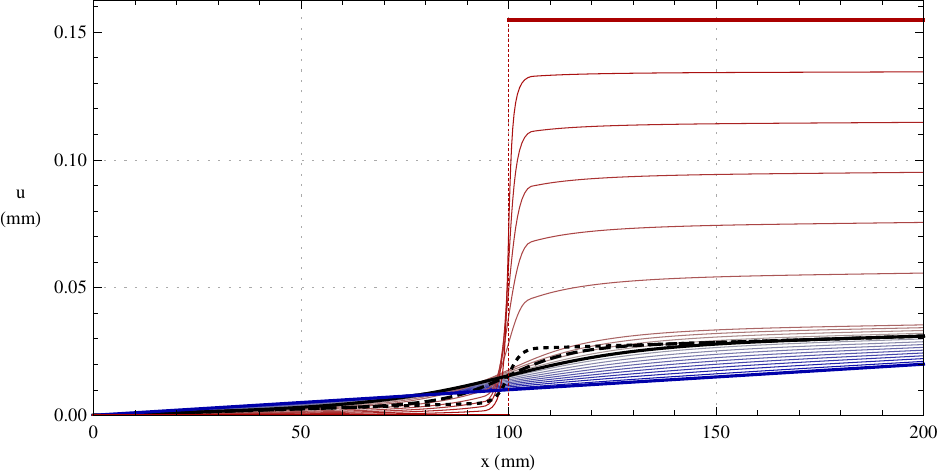}
    \put(12,45){\damBar}
    \put(11,36){\lLeg[0.5]}
  \end{overpic}
  \caption{Displacement profiles at given maximum phase-field states for the bilinear model~\eqref{mod_B} with $\sintL = \SI{10}{mm}$. The displacement profiles at $\sdammax = 0.5$ for smaller internal lengths are also highlighted.}
  \label{fig_UProfB} 
\end{figure}

To define the bilinear softening cohesive law of \Fig{fig_E_L2}, four parameters are necessary, as, for instance, the tensile strength $\scristr$, the fracture toughness $\Gc$, the ratio between the initial critical stress and the intermediate junction stress, and the ratio between the total fracture toughness and the fracture toughness associated with the first linear behaviour. Accordingly we assign
\begin{equation}
  \scristr, \qquad \Gc, \qquad \beta \: \bra{= \dfrac{\scristr}{\tilde\sstr}}, \qquad \gamma \: \bra{= \dfrac{\Gc}{\Gcl}} \,,
  \label{def_L2par} 
\end{equation} 
and deduce the rest of the parameter as
\begin{equation}
\begin{split}
  \Gcl = \gamma \, \Gc, \qquad 
  \tilde\sstr = \beta\,\scristr, \qquad
  \tilde\sope = (1-\beta)\dfrac{2\,\Gcl}{\scristr},  \qquad
  \scriope = 2 \dfrac{(\Gc - \Gcl(1-\beta))}{\beta\scristr}, 
  \\
  k_1 = \dfrac{\scristr^2}{2\,\Gcl}, \qquad
  k_2= \dfrac{\beta^2\scristr^2}{2\bra{\Gc-\Gcl(1-\beta^2)}}, \qquad
  \GcII = \Gc - \Gcl(1-\beta^2) \,.
\end{split}
\end{equation} 
The bilinear traction-separation law reads
\begin{equation}
  \sstr(\sope) =
  \begin{cases}
    \scristr - k_1 \sope, \quad &\text{for $\sope \leq \tilde\sope$} \\
    \tilde\sstr - k_2\bra{\sope - \tilde\sope}, \quad &\text{for $\tilde\sope < \sope \leq\scriope$} \\
    0, \quad &\text{for $\sope > \scriope$}
  \end{cases}.
  \label{eq_TSL_B} 
\end{equation} 
For the third and fourth parameter in \eqref{def_L2par} we respectively assume $\beta=0.3$ and $\gamma=2.5$, this last being a common value adopted for concrete,~\cite{Bazant2002}.
According to \eqref{def_gprime}, the bilinear traction-separation law \eqref{eq_TSL_B} within the mathematical framework becomes
 \begin{equation}
    g'(s)=
    \begin{cases}
      1 - \kscal \, s, \quad &\text{for $s \leq \tilde{s}$} \\
      1 - \kscal \tilde{s} \bra{1-\eta} - \eta \kscal\,s, \quad &\text{for $\tilde{s} < s \leq \bar{s}$} \\
    0, \quad &\text{for $\sope > \scriope$}
  \end{cases}.
  \label{eq_TSLM_B} 
\end{equation}
with $\tilde{s}=(1-\beta/\kscal)$, $\bar{s} = \bra{1 - \kscal \tilde{s} \bra{1-\eta}}/(\eta \kscal)$ and $\eta = k_1/k_2$.
In \cite[Sec. 3.3.3]{Alessi2025b}, an example of material functions $\cbra{\hatf,\PotDann}$ for \eqref{eq_TSLM_B} has been derived 
by fixing $\hatf$, as \eqref{mod_L12}, and deducing, by means of Theorem~3.1 in~\cite{Alessi2025b}, $\PotDann$. For the present case, condition \eqref{con_crack2} provides as scaling factor $\kscal \simeq 0.91/k$.
%
%
The corresponding engineering material functions~$\lDam$ and~$\fDam$ for the considered bilinear model read:
\begin{equation}
  \begin{cases}
    \lDam(\sdam)= (1-\sdam)^2 \\[2ex]
    \fDam(\sdam)= \dfrac{1}{4 k^2} 
    \begin{cases}
      \dfrac{2\sdam-\sdam^2}{\pi^2 \kscal^2}, &\text{ for $0 \leq \sdam < 1-\beta$} 
      \\[1ex]
      \dfrac{
          \bra{ \sqrt{2\sdam-\sdam^2} + \bra{\eta - 1}\sqrt{2\sdam-\sdam^2 + \beta^2 - 1} }^2
        }{\pi^2 \kscal^2}, &\text{ for $1-\beta \leq \sdam \leq 1$} 
    \end{cases}
  \end{cases}.
  \tag{\texttt{B}}
  \label{mod_B} 
\end{equation} 
The degradation function $\fDeg$ and the local material dissipation function $\fDam$ trends for the bilinear model are represented in \Fig{fig_gwB_func}.
The surface fracture energy and the traction-separation law as a function of the displacement jump for the considered bilinear cohesive model are represented in 
Fig.s~\ref{GBPlot} and~\ref{tslawBPlot1}, respectively. 
The global response of \Fig{fig_UFB} clearly reflects the bilinear cohesive fracture behaviour. The initial linear elastic response is indeed followed by a bilinear softening damage response. Also for the bilinear model, the global response does not change with the internal length. The internal length, if sufficiently small compared to the bar length, only affects the phase-field localization size, as evident from the phase-field and displacement field profiles evolutions represented, respectively, in \Fig{fig_DamProfB} and \Fig{fig_UProfB}, where the same three internal lengths~\eqref{def_lL} of the linear model have been considered.
The trend of the half-width phase-field localization support~\eqref{eq_D} with respect to the fracture evolution is highlighted in \Fig{fig_DB}.
To the best of the author's knowledge, no phase-field model capable of capturing a bilinear softening cohesive response has yet been presented in the literature.

%
%

\FloatBarrier
\subsection{Exponential softening cohesive law}
\label{sec_E}

The exponential or Barenblatt like cohesive law, widely adopted in applications and also considered in \cite{Wu2017}, is given by the following traction-separation law
\begin{equation}
  \sstr(\sope) =   \scristr \, \exp\bra{-\frac{\scristr}{\Gc}\sope}.
  \label{eq_TSL_E} 
\end{equation} 
Its behaviour is fully determined by two parameters, the critical stress and the fracture toughness, and is characterized by an infinite ultimate crack opening~\eqref{def_scriope}, namely $\scriope = +\infty$.
According to \eqref{def_gprime}, the linear traction-separation law within the mathematical framework \eqref{eq_TSL_E} becomes
\begin{equation}
  g'(s)= \exp\bra{-\kscal s} \,.
  \label{eq_TSLM_E} 
\end{equation}

In \cite[Sec. 3.3.6]{Alessi2025b}, an example of material functions $\cbra{\hatf,\PotDann}$ for \eqref{eq_TSLM_E} has been derived by fixing $\hatf$ and deducing, by means of Theorem~3.1 in~\cite{Alessi2025b}, $\PotDann$, with condition \eqref{con_crack2} providing as scaling factor $\kscal = 1/k$.
In fact, a slightly modified version of \eqref{eq_TSLM_E} was considered in \cite[Sec. 3.3.6]{Alessi2025b} to ensure that all the assumptions of \cite[Theorem~1.1]{Alessi2025a} were satisfied, specifically the continuity of the pair $\cbra{\hatf,\PotDann}$. This modified traction-separation law introduces a small parameter that influences the fracture behavior only at very large crack opening displacements.
Neglecting this modification, which affects only the asymptotic behavior, the corresponding engineering material functions $\lDam$ and $\fDam$ for the exponential model \eqref{eq_TSL_E} simplify to:
\begin{equation}
  \begin{cases}
    \lDam(\sdam)= (1-\sdam)^2 \\[2ex]
    \fDam(\sdam)= \dfrac{\bra{\cosh^{-1}\bra{\frac{1}{1-\alpha }}}^2}{4\,\pi^2}
  \end{cases} \,.
  \tag{\texttt{E}}
  \label{mod_E} 
\end{equation} 
The following analysis indicates that the impact on practical applications of the modified version of the model is negligible and that the exponential cohesive fracture response \eqref{eq_TSL_E} based on the material functions \eqref{mod_E} is well captured. 

Unlike in \cite{Wu2017}, where the exponential cohesive law was derived by approximation with fitted polynomial functions, \eqref{mod_E} provides an analytical closed-form expression. Being $\lim_{\sdam \to 1}\fDam(\sdam) = +\infty$, cohesive forces never vanish across the crack.

\begin{figure}[h]
\centering
\begin{subfigure}[b]{0.4\linewidth}
  \includegraphics[width=0.8\linewidth]{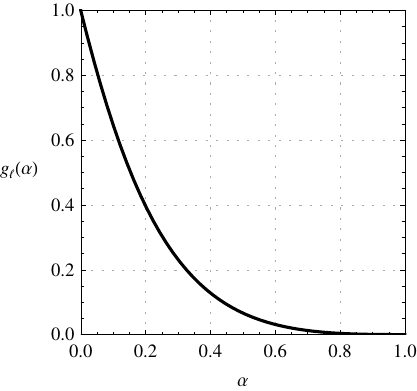}
  \caption{}
  \label{gEPlot} 
\end{subfigure}
\begin{subfigure}[b]{0.4\linewidth}
  \includegraphics[width=0.8\linewidth]{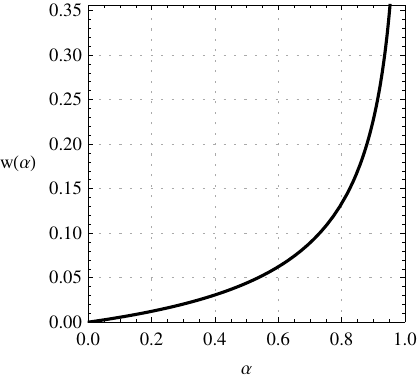}
  \caption{}
  \label{wEPlot} 
\end{subfigure}
  \caption{(\subref{gEPlot}) Degradation and (\subref{wEPlot}) local dissipation functions for the exponential cohesive fracture model~\eqref{mod_E}.}
  \label{fig_gwE_func} 
\end{figure}

The degradation function $\fDeg$ and the local material dissipation function $\fDam$ trends for the exponential model~\eqref{mod_E} are represented in \Fig{fig_gwE_func}.

\begin{figure}[h]
\centering 
\begin{subfigure}[b]{0.31\linewidth}
  \includegraphics[width=\linewidth]{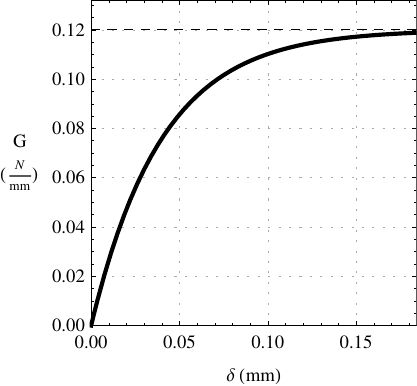}
  \caption{}
  \label{GEPlot} 
\end{subfigure}
\quad
\begin{subfigure}[b]{0.31\linewidth}
  \begin{overpic}[width=\linewidth]{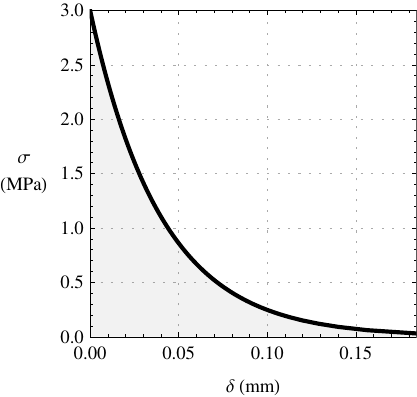}
     \put(29,25){$\Gc$}
  \end{overpic}
  \caption{}
  \label{tslawEPlot1} 
\end{subfigure}
\\[2ex]
\begin{subfigure}[b]{0.31\linewidth}
  \includegraphics[width=\linewidth]{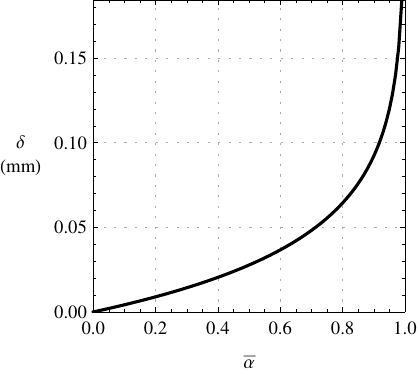}
  \caption{}
  \label{tslawEPlot2} 
\end{subfigure}
\quad
\begin{subfigure}[b]{0.31\linewidth}
  \includegraphics[width=\linewidth]{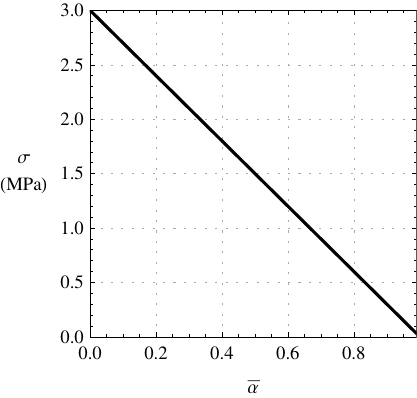}
  \caption{}
  \label{tslawEPlot3} 
\end{subfigure}
  \caption{(\subref{GEPlot}) Surface fracture energy and (\subref{tslawEPlot1}) traction-separation law of the exponential model~\eqref{mod_E} as a function of the displacement jump, obtained using the parametrization \eqref{def_param} of the displacement jump (\subref{tslawEPlot2}) and the critical stress~(\subref{tslawEPlot3}) with respect to the maximum phase-field value.}
  \label{fig_tslawE} 
\end{figure}

Fig.s~\ref{GEPlot} and \ref{tslawEPlot1} respectively represent the surface fracture energy and the traction-separation law as a function of the displacement jump for the considered exponential cohesive model.

\begin{figure}[h]
\centering
  \begin{overpic}[scale=0.8]{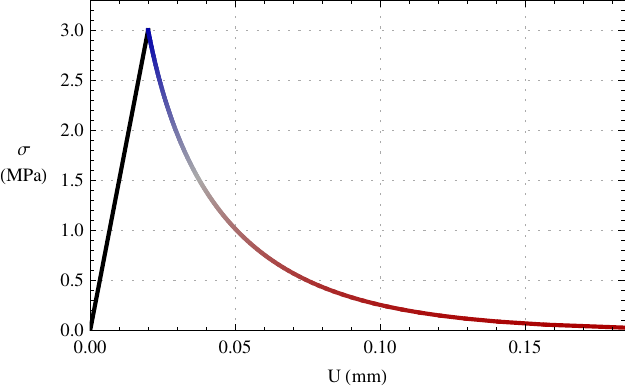}
    \put(76,56){\damBar[$\sdammax$]}
  \end{overpic}
  \caption{Global response for the exponential model~\eqref{mod_E}, independent of the internal length. The initial elastic response is shown in black. The dependence on the maximum phase-field state is also highlighted.}
  \label{fig_UFE} 
\end{figure}

\begin{figure}[h]
\centering
  \begin{overpic}[width=0.8\linewidth]{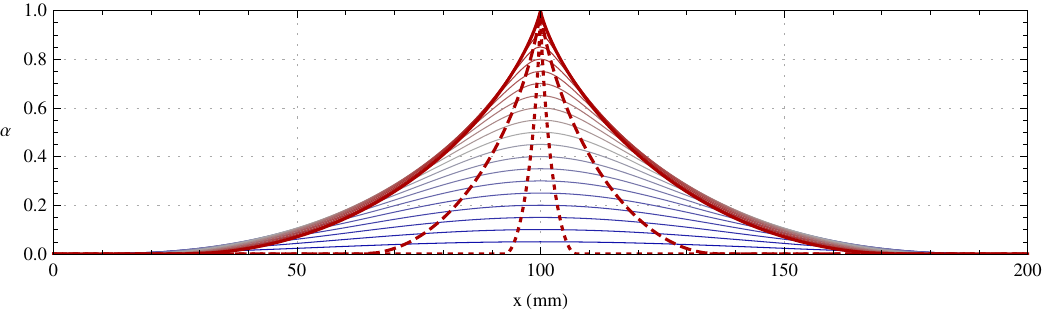}
    \put(81,24){\damBar}
    \put(76,14){\lLeg}
  \end{overpic}
  \caption{Phase-field profiles at given maximum phase-field states for the exponential model~\eqref{mod_E}. The ultimate phase-field profiles associated with the three internal lengths~\eqref{def_lL} are also highlighted.
  }
  \label{fig_DamProfE} 
\end{figure}

\begin{figure}[h]
\centering
\footnotesize
  \includegraphics[scale=0.9]{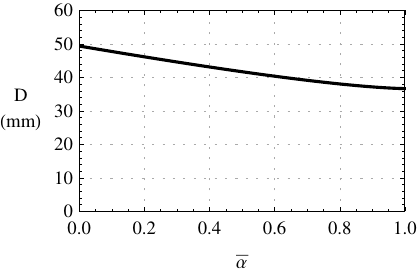}
  \caption{Half-width of the support of the phase-field profile for the exponential model~\eqref{mod_E} with respect to the fracture evolution ($\sintL=\SI{5}{mm}$).}
  \label{fig_DE} 
\end{figure}

\begin{figure}[h]
\centering
  \begin{overpic}[width=0.8\linewidth]{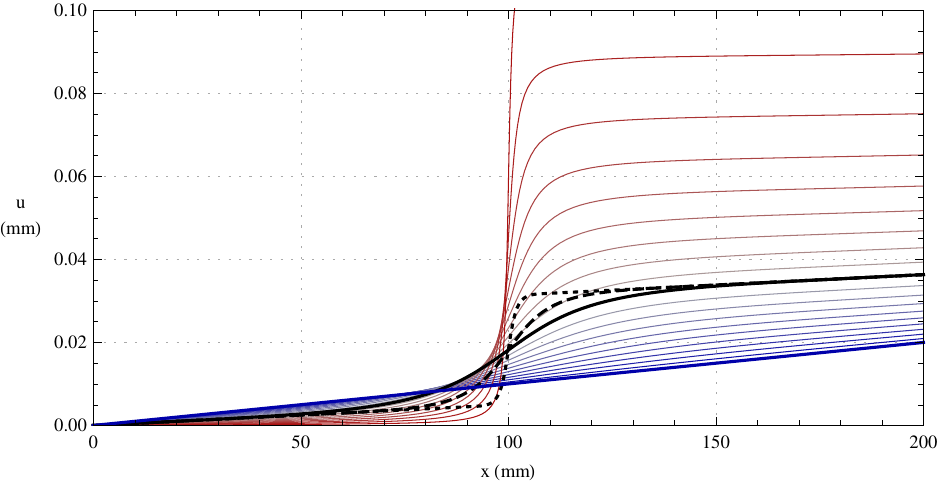}
    \put(12,45){\damBar}
    \put(11,35){\lLeg[0.5]}
  \end{overpic}
  \caption{Displacement profiles at given maximum phase-field states for the exponential model~\eqref{mod_E} with $\sintL = \SI{5}{mm}$. The displacement profiles at $\sdammax = 0.5$ for smaller internal lengths are also highlighted.}
  \label{fig_UProfE} 
\end{figure}

The global response of \Fig{fig_UFE} clearly reflects the exponential cohesive fracture behaviour. The initial linear elastic response is indeed followed by a concave softening damage response approaching asymptotically zero. Also for the exponential model, the global response does not change with the internal length. The internal length, if sufficiently small compared to the bar length, only affects the phase-field localization size, as evident from the phase-field and displacement field profiles evolutions represented, respectively, in Fig.s \ref{fig_DamProfE} and \ref{fig_UProfE}, where the three internal lengths~\eqref{def_lL} have been considered. 
The trend of the half-width phase-field localization support~\eqref{eq_D} with respect to the fracture evolution is highlighted in \Fig{fig_DE}.

\FloatBarrier
\subsection{Hyperbolic softening cohesive law}
\label{sec_H}

In this section we consider two hyperbolic cohesive fracture models, both defined through only two parameters, which mainly differs for their limit behaviour. 

\paragraph{Linear hyperbolic function}
The engineering traction-separation law of the first hyperbolic model is given by 
\begin{equation}
  \sstr(\sope) = 
  \begin{cases}
    \scristr \bra{\dfrac{2}{1+\sope/\scriope}-1}, \quad &\text{for $\sope \leq \scriope$} \\
    0, \quad &\text{for $\sope > \scriope$} 
  \end{cases} \,,
  \label{eq_TSL_H1} 
\end{equation} 
with $ \scriope = \Gc/\bra{\scristr(2\ln 2 -1)}$ being the ultimate crack opening.
According to \eqref{def_gprime}, the traction-separation law within the mathematical framework for \eqref{eq_TSL_H1} becomes
\begin{equation}
  g'(s)=
  \begin{cases}
      \dfrac{2}{1+\kscal\,s}-1, \; & \text{for $s \leq \frac{1}{\kscal}$} \\
      0, \; & \text{for $s > \frac{1}{\kscal}$}
  \end{cases} \,.
\label{eq_TSLM_H1} 
\end{equation}
In \cite[Sec. 3.3.4]{Alessi2025b}, an example of material functions $\cbra{\hatf,\PotDann}$ for \eqref{eq_TSLM_H1} has been derived by fixing $\hatf$ and deducing, by means of Theorem~3.1 in~\cite{Alessi2025b}, $\PotDann$, with condition \eqref{con_crack2} providing as scaling factor $\kscal \simeq 0.386/k$.
The corresponding engineering material functions~$\lDam$ and~$\fDam$ for the considered hyperbolic model read:
\begin{equation}
  \begin{cases}
    \lDam(\sdam)= (1-\sdam)^2 \\[2ex]
    \fDam(\sdam) = \dfrac{0.170 \, \bra{ 2\sdam - \sdam^2 + 2 (1 -\sdam)^2 \log\bra{1-\sdam}}^2}{(2 - \sdam)^3 \sdam^3}
  \end{cases}.
  \tag{\texttt{H1}}
  \label{mod_H1} 
\end{equation}

\paragraph{Quadratic hyperbolic function}
The engineering traction-separation law of the second hyperbolic model is given by
\begin{equation}
  \sstr(\sope) = 
     \dfrac{\scristr}{\bra{1+  \sope/\hat\sope}^2} \,,
  \label{eq_TSL_H2} 
\end{equation} 
with $\hat\sope = \Gc/\scristr$.
According to \eqref{def_gprime}, the traction-separation law within the mathematical framework for \eqref{eq_TSL_H2}
becomes
\begin{equation}
  g'(s)= \dfrac{1}{\bra{1+\kscal\,s}^2} \,.
  \label{eq_TSLM_H2} 
\end{equation}
In \cite[Sec. 3.3.5]{Alessi2025b}, an example of material functions $\cbra{\hatf,\PotDann}$ for \eqref{eq_TSLM_H2} has been derived by fixing $\hatf$ and deducing, by means of Theorem~3.1 in~\cite{Alessi2025b}, $\PotDann$, with condition \eqref{con_crack2} providing as scaling factor $\kscal \simeq 1/k$.
The corresponding engineering material functions~$\lDam$ and~$\fDam$ for the considered hyperbolic model read:
\begin{equation}
  \begin{cases}
    \lDam(\sdam)= (1-\sdam)^2 \\[2ex]
    \fDam(\sdam)= \dfrac{
      \bra{
      _2F_1\bra{-1/4,1,1/2,2\sdam - \sdam^2} - 
      (1-\sdam)^2\, _2F_1\bra{3/4,1,1/2,2\sdam - \sdam^2}
      }^2
        }{
      (4\pi)^2 \bra{2\sdam - \sdam^2}}
  \end{cases} \,,
  \tag{\texttt{H2}}
  \label{mod_H2} 
\end{equation} 
with $_2F_1$ being the hypergeometric function.

The main difference between the two hyperbolic models is that model~\eqref{eq_TSL_H1} owns a finite ultimate crack opening whereas \eqref{eq_TSL_H2} does not. The second model has also been considered in~\cite{Wu2017}, although with material functions obtained by a fitting procedure. None of the two hyperbolic models have been instead considered in \cite{Feng2021}.

\begin{figure}[h]
\centering
\begin{subfigure}[b]{0.4\linewidth}
  \begin{overpic}[width=0.8\linewidth]{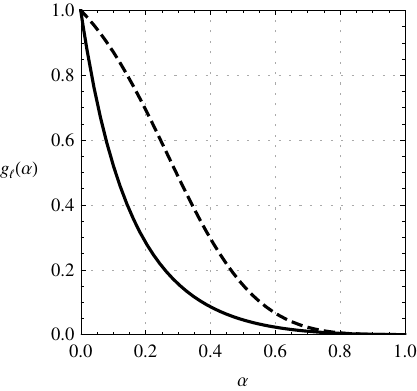}
    \put(59,75){\HLeg}
  \end{overpic}
  \caption{}
  \label{gHPlot} 
\end{subfigure}
\begin{subfigure}[b]{0.4\linewidth}
  \begin{overpic}[width=0.8\linewidth]{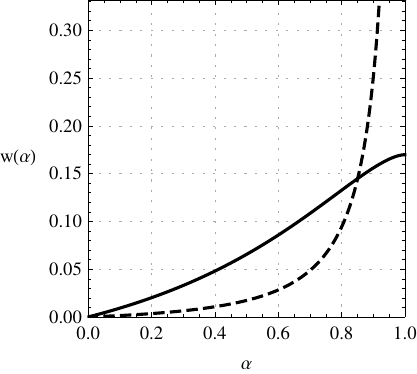}
    \put(22,75){\HLeg}
  \end{overpic}
  \caption{}
  \label{wHPlot} 
\end{subfigure}
  \caption{(\subref{gHPlot}) Degradation and (\subref{wHPlot}) local dissipation functions for the hyperbolic cohesive fracture models.}
  \label{fig_gwH_func} 
\end{figure}

\begin{figure}[h]
\centering 
\begin{subfigure}[b]{0.31\linewidth}
  \begin{overpic}[width=\linewidth]{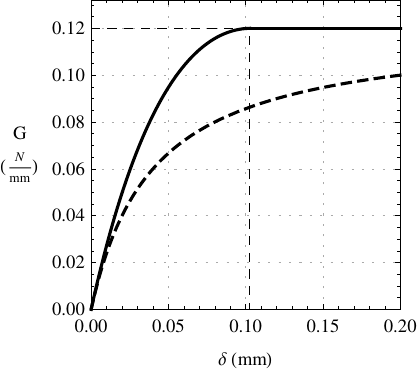}
    \put(60,22){$\HLeg$}
  \end{overpic}
  \caption{}
  \label{GHPlot} 
\end{subfigure}
\quad
\begin{subfigure}[b]{0.31\linewidth}
  \begin{overpic}[width=\linewidth]{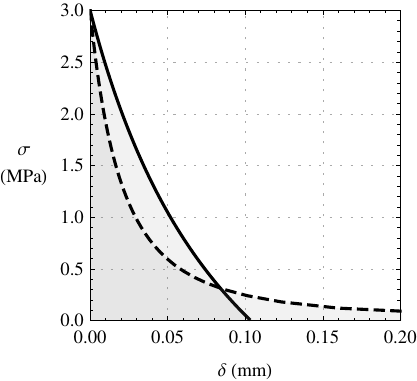}
     \linethickness{.1pt}  
     \put(35,45){\line(1,1){10}}  
     \put(25,19){$\Gc$ {\footnotesize(\texttt{H2})}}
     \put(46,56){$\Gc$ {\footnotesize(\texttt{H1})}}
     \put(60,75){$\HLeg$}
  \end{overpic}
  \caption{}
  \label{tslawHPlot1} 
\end{subfigure}
\\[2ex]
\begin{subfigure}[b]{0.31\linewidth}
  \begin{overpic}[width=\linewidth]{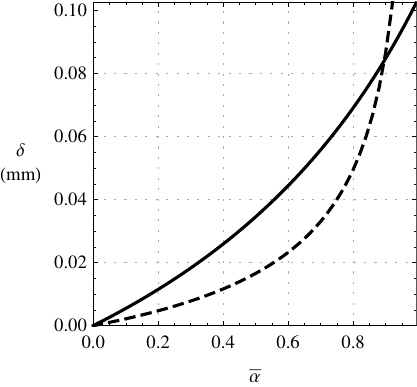}
    \put(24,78){$\HLeg$}
  \end{overpic}
  \caption{}
  \label{tslawHPlot2} 
\end{subfigure}
\quad
\begin{subfigure}[b]{0.31\linewidth}
  \begin{overpic}[width=\linewidth]{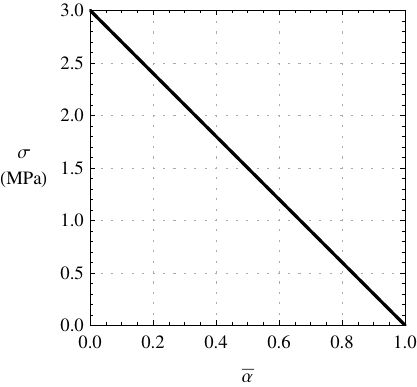}
    \put(60,75){$\HLeg$}
  \end{overpic}
  \caption{}
  \label{tslawHPlot3} 
\end{subfigure}
  \caption{(\subref{GHPlot}) Surface fracture energy and (\subref{tslawHPlot1}) traction-separation law as a function of the displacement jump, obtained using the parametrisation \eqref{def_param} of the displacement jump (\subref{tslawHPlot2}) and the critical stress (\subref{tslawHPlot3}) with respect to the maximum phase-field value for the considered hyperbolic models \eqref{mod_H1} and \eqref{mod_H2}.}
  \label{fig_GH} 
  \label{fig_tslawH} 
\end{figure}

\begin{figure}[h]
\centering
  \begin{overpic}[scale=0.8]{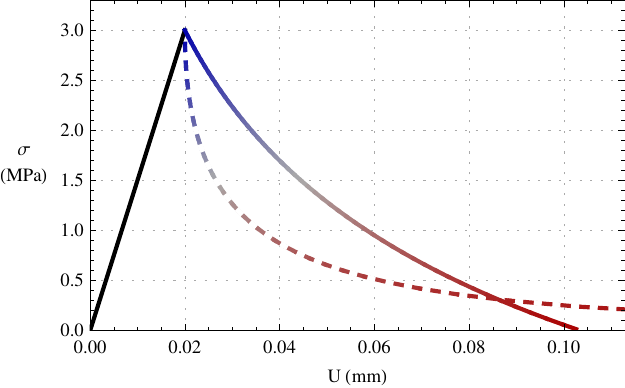}
    \put(77,53){\damBar[$\sdammax$]}
    \put(78,42){$\HLeg$}
  \end{overpic}
    \caption{Global response of the bar for the hyperbolic models (\eqref{mod_H1} and \eqref{mod_H2}), independent of the internal length. The initial elastic response is shown in black. The dependence on the maximum phase-field state is also highlighted.}
  \label{fig_UFH} 
\end{figure}

In \Fig{fig_gwH_func}, the degradation function $\fDeg$ and the local material dissipation function $\fDam$ trends for the two hyperbolic models are depicted.
Fig.s~\ref{GHPlot} and~\ref{tslawHPlot1} respectively represent the surface fracture energy and the traction-separation law as a function of the displacement jump for all the considered hyperbolic models. 
Their global responses, which again are independent from the internal length, are compared in \Fig{fig_UFH}.

To highlight the influence of the~$\sintL$, a set of three different values for the internal length for the two hyperbolic models have been considered:
\begin{subequations} \label{def_lH} 
  \begin{align} 
    &\text{\eqref{mod_H1}: $\sintL = \cbra{10,5,1} \quad \text{(mm)}$} \,, \label{def_lH1} \\
    &\text{\eqref{mod_H2}: $\sintL = \cbra{5,2.5,0.5} \quad \text{(mm)}$} \,. \label{def_lH2} 
  \end{align}           
\end{subequations}

\begin{figure}[h]
\centering
\footnotesize
  \begin{overpic}[scale=0.9]{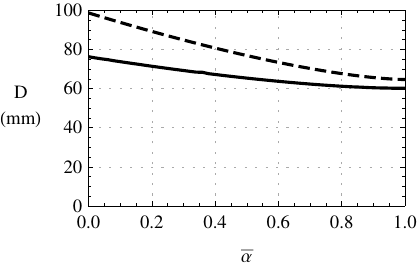}
    \put(68,21){\HLeg}    
  \end{overpic}
  \caption{Half-width of the support of the phase-field profiles for the hyperbolic models with respect to the fracture evolution (\eqref{mod_H1}: $\sintL=\SI{10}{mm}$ and \eqref{mod_H2}: $\sintL=\SI{5}{mm}$).}
  \label{fig_DH} 
\end{figure}

\begin{figure}[h]
\centering
\begin{subfigure}{\linewidth}
  \centering
  \begin{overpic}[width=0.8\linewidth]{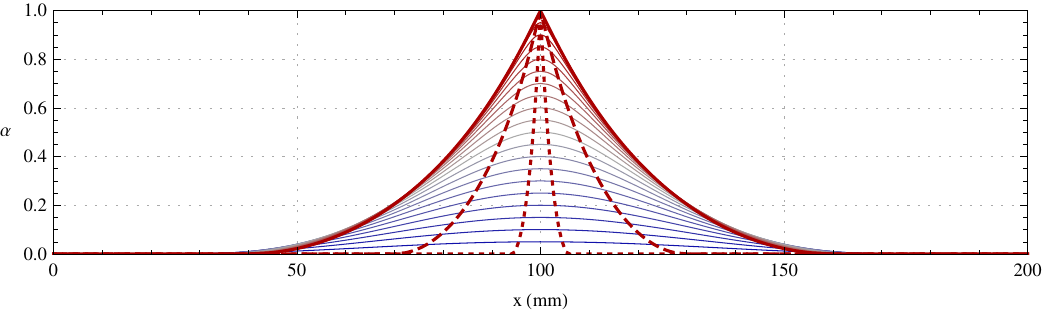}
    \put(82,23){\damBar}
    \put(77,13){\lLeg}
  \end{overpic}
  \caption{}
  \label{fig_DamProfH1} 
\end{subfigure}
\\[2ex]
\begin{subfigure}{\linewidth}
  \centering
  \begin{overpic}[width=0.8\linewidth]{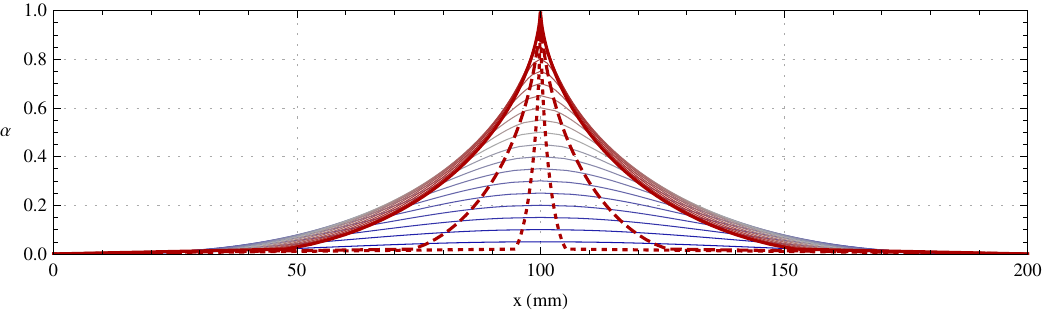}
    \put(82,23){\damBar}
    \put(77,13){\llLeg}
  \end{overpic}
  \caption{}
  \label{fig_DamProfH2}
\end{subfigure}
  \caption{Phase-field profiles at given maximum phase-field states for the hyperbolic models:~(\subref{fig_DamProfH1}) for~\eqref{mod_H1} and~(\subref{fig_DamProfH2}) for~\eqref{mod_H2}. The ultimate phase-field profiles associated with different internal lengths are also highlighted.}
  \label{fig_DamProfH} 
\end{figure}

\Fig{fig_DamProfH} reports the phase-field profiles evolutions, obtained from \eqref{eq_damprof}, together with the ultimate phase-field profiles associated with the internal lengths~\eqref{def_lH}.
The trends of the half-width localization support~\eqref{eq_D} with respect to the fracture evolution are highlighted in \Fig{fig_DH}.

\begin{figure}[h]
\centering
\begin{subfigure}{\linewidth}
  \centering
  \begin{overpic}[width=0.8\linewidth]{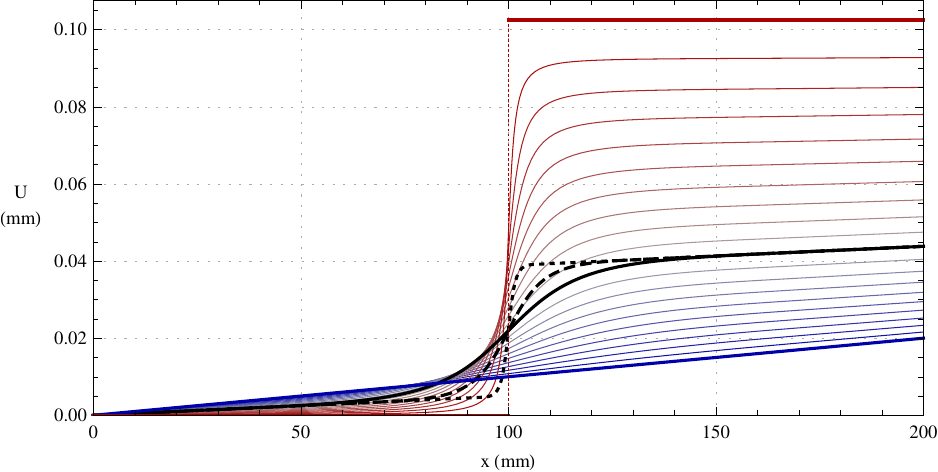}
    \put(11,45){\damBar}
    \put(11,35){\lLeg[0.5]}
  \end{overpic}
  \caption{}
  \label{fig_UProfH1} 
\end{subfigure}
\\[2ex]
\begin{subfigure}{\linewidth}
  \centering
  \begin{overpic}[width=0.8\linewidth]{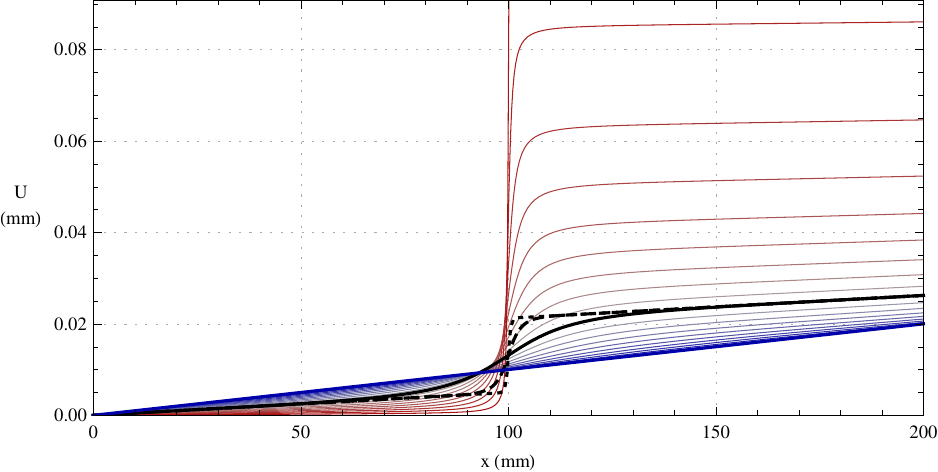}
    \put(11,45){\damBar}
    \put(11,35){\llLeg[0.5]}
  \end{overpic}
  \caption{}
  \label{fig_UProfH2}
\end{subfigure}
  \caption{
    Displacement profiles at given maximum phase-field states for the hyperbolic models:~(\subref{fig_UProfH1}) for~\eqref{mod_H1} with~$\sintL = \SI{10}{mm}$   and~(\subref{fig_UProfH2}) for~\eqref{mod_H2} with $\sintL = \SI{5}{mm}$. The displacement profiles at $\sdammax = 0.5$ for smaller internal lengths are also highlighted.
  }
  \label{fig_UProfH} 
\end{figure}

The evolution of the displacement profiles for \eqref{mod_H1} and \eqref{mod_H2} are reported in \Fig{fig_UProfH}. Therein, the convergence of the smooth displacement profiles towards profiles with a discontinuity (the crack) is noticeable as the internal length decreases.

\clearpage
\subsection{Dugdale's cohesive law} 
\label{sec_D} 


\newcommand{\DMLeg}{
\colorbox{white}{
\begin{minipage}{10mm}
\scriptsize
  \includegraphics[trim=2mm 1.8mm 0mm 0mm,clip,scale=0.8]{L1}\:D1 \\
  \includegraphics[trim=2mm 1.8mm 0mm 0mm,clip,scale=0.8]{L2}\:D2
\end{minipage}
}}

\newcommand{\DellLeg}{
\colorbox{white}{
\begin{minipage}{21mm}
\scriptsize
  \scalebox{1}[1.000]{\includegraphics[trim=2mm 1.8mm 0mm 0mm,clip,scale={0.8}]{L1}}\:$\ell=\SI{10}{mm}$ \\
  \scalebox{1}[0.666]{\includegraphics[trim=2mm 1.8mm 0mm 0mm,clip,scale={0.8}]{L1}}\:$\ell=\SI{5}{mm}$ \\
  \scalebox{1}[0.333]{\includegraphics[trim=2mm 1.8mm 0mm 0mm,clip,scale={0.8}]{L1}}\:$\ell=\SI{1}{mm}$
\end{minipage}
}}


As the linear model, Dugdale cohesive law is completely defined once two of the following three parameters are defined: the tensile strength $\scristr$, the fracture toughness $\Gc$ and the ultimate crack opening~$\scriope$ since $\Gc = \scriope\,\scristr$, \Fig{fig_E_D}. The Dugdale traction-separation law reads
\begin{equation}
  \sstr(\sope) =
  \begin{cases}
    \scristr , \quad &\text{for $\sope \leq \scriope$} \\
    0, \quad &\text{for $\sope > \scriope$} 
  \end{cases} \,,
  \label{eq_TSL_D} 
\end{equation} 
where the ultimate crack opening is given by $\scriope = \Gc/\scristr = \SI{0.04}{\mm}$ with the constitutive parameters of \Tab{tab_par}.
According to \eqref{def_gprime}, the linear traction-separation law within the mathematical framework \eqref{eq_TSL_D} becomes
 \begin{equation}
    g'(s)=
    \begin{cases}
        1, \; & \text{for $s \leq \kscal$} \\
        0, \; & \text{for $s > \kscal$}
    \end{cases} \,.
  \label{eq_TSLM_D} 
\end{equation}
Taking advantage of the results of Theorem~2.18 in \cite[Sec. 3.3.1]{Alessi2025b}, to which Dugdale model applies,
two different models have been derived, both describing the same Dugdale cohesive behaviour and sharing the same scaling factor $\kscal=k$. 
Compared to Theorems~3.1 and~3.2 in \cite{Alessi2025b}, Theorem~2.18 provides a slightly less constructive procedure for the deduction of the material functions $\hatf$ and $\PotDann$ of the phase-field model.
The engineering material functions $\lDam$ and $\fDam$ for the present two models read:


\begin{equation}
\begin{cases}
  \lDam^{-1}(t)= 1-\bra{1-\dfrac{2}{\pi}\bra{\arcsin(\sqrt{t}) - \sqrt{t-t^2}} }^{2/3} \\
  \fDam(\sdam)= \dfrac{9}{64}\sdam
\end{cases},
\tag{\texttt{D}$_1$}
\label{mod_D1} 
\end{equation} 
and
\begin{equation}
\begin{cases}
  \lDam^{-1}(t) = 1-\bra{1-\dfrac{2}{\pi}\bra{\arcsin(\sqrt{t}) - \sqrt{t-t^2}} }^{1/2} \\
  \fDam(\sdam) = \dfrac{\sdam^2}{4}
\end{cases}.
\tag{\texttt{D}$_2$}
\label{mod_D2} 
\end{equation}

Both models were derived by prescribing $\PotDann$, linear in \eqref{mod_D1} and quadratic in \eqref{mod_D2}, and assuming $\Psi$, from which (the inverse of) $\hatf$ was obtained.
The function $\Psi$ in \cite[$(2.1)$]{Alessi2025b} was chosen so that the auxiliary function (2.30) in \cite{Alessi2025b} corresponds to $\Phi(x) = 2\,\kscal\,\hatf^{1/2}(x)$.

\begin{figure}[h]
\centering
\begin{subfigure}[b]{0.4\linewidth}
  \begin{overpic}[width=0.8\linewidth]{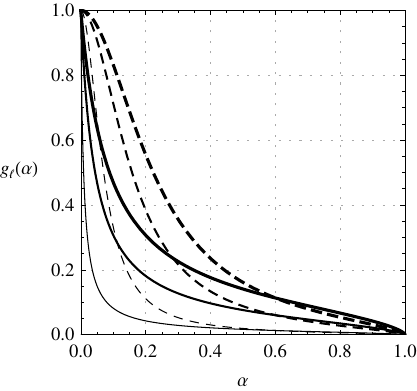}
    \put(53,80){\DMLeg}
    \put(44,62){\DellLeg}
  \end{overpic}
  \caption{}
  \label{gDPlot} 
\end{subfigure}
\begin{subfigure}[b]{0.4\linewidth}
  \begin{overpic}[width=0.8\linewidth]{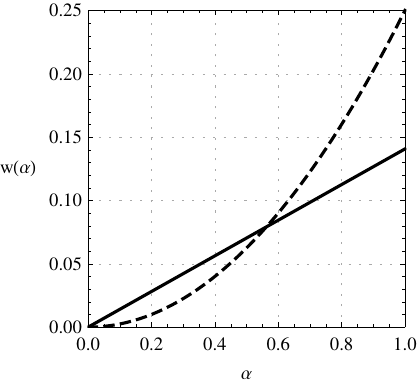}
    \put(22,75){\DMLeg}
  \end{overpic}
  \caption{}
  \label{wDPlot} 
\end{subfigure}
  \caption{(\subref{gDPlot}) Degradation and (\subref{wDPlot}) local dissipation functions for the Dugdale cohesive fracture models~\eqref{mod_D1} and~\eqref{mod_D2}, for different internal lengths.}
  \label{fig_gwD_func} 
\end{figure}

As for the linear models, the presence or absence of a linear term in the phase-field local dissipation function is not a necessary condition for an explicit elastic response.
The function $\lDam(\sdam)$ for models \eqref{mod_D1} and \eqref{mod_D2} does not admit an explicit analytical inverse. As a consequence, an interpolation function to reconstruct $\lDam(\sdam)$ is hereafter adopted, similar to what have been done for \eqref{mod_L1} and \eqref{mod_L2}.
The degradation function $\fDeg$ and the local material dissipation function $\fDam$ trends for the two models are depicted in \Fig{fig_gwD_func}. Therein the influence of the internal length is also highlighted.

\begin{figure}[h]
\centering
\begin{subfigure}[b]{0.31\linewidth}
  \begin{overpic}[width=\linewidth]{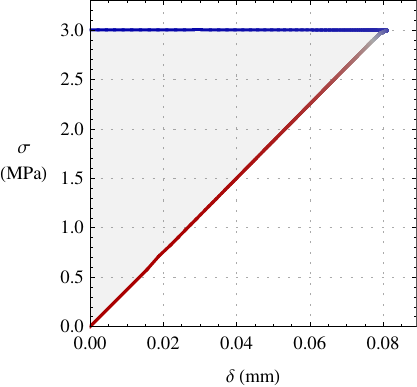}
        \put(35,60){$\Gc$}
        \put(63,40){\damBar}
        \put(64,23){\DMLeg}
  \end{overpic}
  \caption{}
  \label{tslawDPlot1} 
\end{subfigure}
\quad
\begin{subfigure}[b]{0.31\linewidth}
  \begin{overpic}[width=\linewidth]{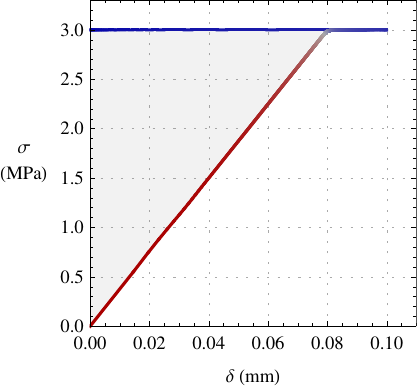}
        \put(35,60){$\Gc$}
        \put(63,40){\damBar}
        \put(64,23){\DMLeg}
  \end{overpic}
  \caption{}
  \label{tslawDPlot2} 
\end{subfigure}
\\[2ex]
\begin{subfigure}[b]{0.31\linewidth}
  \begin{overpic}[width=\linewidth]{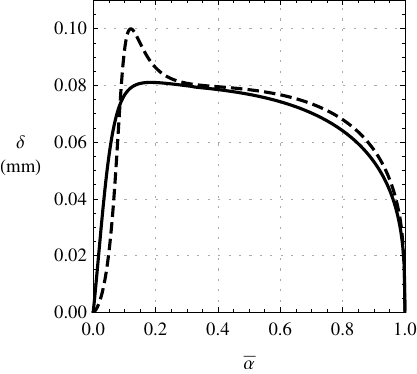}
    \put(40,30){\DMLeg}
  \end{overpic}
  \caption{}
  \label{tslawDPlot3} 
\end{subfigure}
\quad
\begin{subfigure}[b]{0.31\linewidth}
  \begin{overpic}[width=\linewidth]{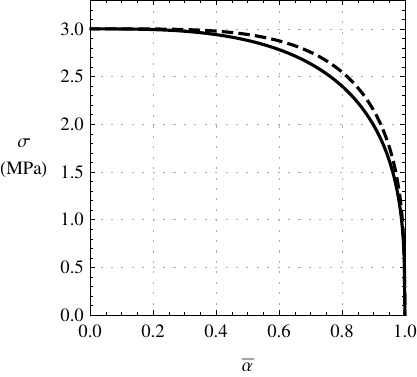}
    \put(40,30){\DMLeg}
  \end{overpic}
  \caption{}
  \label{tslawDPlot4} 
\end{subfigure}
  \caption{Traction-separation laws (\subref{tslawDPlot1}) and (\subref{tslawDPlot2}) for, respectively \eqref{mod_D1} and  \eqref{mod_D2} as a function of the displacement jump, obtained using the parametrisation \eqref{def_param} of the displacement jump (\subref{tslawDPlot3}) and the critical stress (\subref{tslawDPlot4}) with respect to the maximum phase-field value. In (\subref{tslawDPlot1}) and (\subref{tslawDPlot2}) the dependence on the maximum phase-field state is also highlighted.
  }
  \label{fig_GD} 
\end{figure}

Fig.s \ref{tslawDPlot1} and \ref{tslawDPlot2} respectively represent the traction-separation laws as a function of the displacement jump for the two considered Dugdale models. The link between these responses with the maximum value of the phase-field is also provided, deduced from the parametrisation \eqref{def_param} of the displacement jump (\Fig{tslawDPlot3}) and the critical stress (\Fig{tslawDPlot4}).

\begin{figure}[h]
\begin{subfigure}[b]{\linewidth}
  \centering
  \begin{overpic}[scale=0.8]{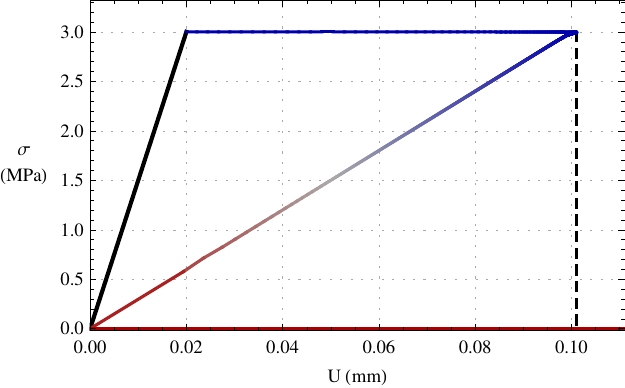}
    \put(70,15){\damBar}
  \end{overpic}
  \caption{}
  \label{UFDPlot1} 
\end{subfigure}
\\[2ex]
\begin{subfigure}[b]{\linewidth}
  \centering
  \begin{overpic}[scale=0.8]{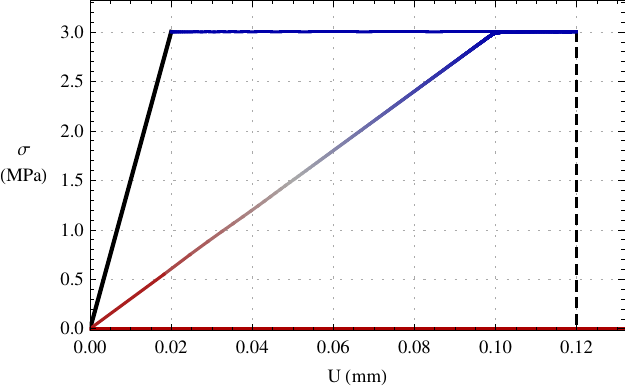}
    \put(70,15){\damBar}
  \end{overpic}
  \caption{}
  \label{UFDPlot2} 
\end{subfigure}
  \caption{Global responses (\subref{UFDPlot1}) and (\subref{UFDPlot2}) for respectively \eqref{mod_D1} and \eqref{mod_D2}, independent of the internal length. The initial elastic response is shown in black. The dependence on the maximum phase-field state is also highlighted. At the critical displacement~$U^*$, a snap back response occur with the critical stress suddenly dropping from the limit stress value to zero.}
  \label{fig_UFD} 
\end{figure}

Differently from the previous cohesive fracture responses, the displacement jump is not monotonically increasing with respect to the maximum phase-field value. Let us identify the critical crack opening, that is the maximum observable crack opening associated with non-vanishing cohesive forces, and the associated critical phase-field value as
\begin{equation}
  \sope^* \coloneqq \max_{\sdam}\sope(\sdam) 
  \qquad \text{and} \qquad
  \sdam^* \coloneqq \argmax_{\sdam}\sope(\sdam)  \,.
\end{equation} 
As a consequence, the energy dissipated within the bar is not equal to the fracture toughness, although \eqref{eq_Gc} continues to hold, but rather given by $G=G(\sope^*)=\Gc\,\sope^*$.
Indeed, the global responses in Fig.s~\ref{UFDPlot1} and~\ref{UFDPlot2} highlights the occurrence of a snap back response. At the critical state
$\cbra{U^*,\sope^*,\sdam^*}\simeq\cbra{\SI{0.101}{mm},\SI{0.081}{mm},0.186}$ for \eqref{mod_D1} and 
$\cbra{U^*,\sope^*,\sdam^*}\simeq\cbra{\SI{0.120}{mm},\SI{0.100}{mm},0.120}$ for \eqref{mod_D2},
the phase-field variable suddenly localizes to its optimal profile for $\sdammax = 1$. Except for the critical displacement value, the two models describe the same global response. The phenomenology is essentially equal to that one observed in the same problem for \texttt{AT}$_1$ and \texttt{AT}$_2$ phase-field brittle fracture models, \cite{Tanne2018}, where an instable brittle fracture occurs at the peak load. Nevertheless, the present phase-field approximations correctly describe the sought cohesive response in terms of crack opening--limit stress evolution and global response.

\begin{figure}[h]
\centering
\footnotesize
  \begin{subfigure}[b]{0.6\linewidth}
    \begin{overpic}[scale=0.9]{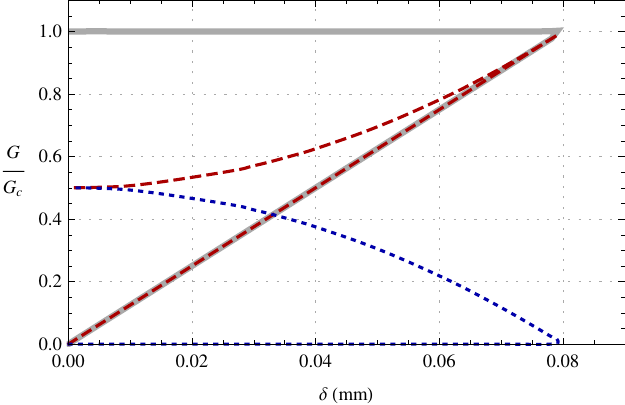}
      \put(23,22){\rotatebox{-57}{\normalsize\color{red!50!black}$\boldsymbol\uparrow$}}
      \put(30,11){\rotatebox{1}{\normalsize\color{blue}$\boldsymbol\rightarrow$}}
      \put(66,31){
      \colorbox{white}{
        \begin{minipage}{25mm}
        \scriptsize
          \includegraphics[trim=0mm 1mm -1mm 0mm,clip,scale=0.8]{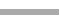}\:$\int_{0}^{L}\gamma(\sdam,\sdam')\dx$\\[2ex]
          \includegraphics[trim=0mm 1mm -1mm 0mm,clip,scale=0.8]{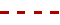}\:$\int_{0}^{L}\fDam(\sdam)/\sintL\dx$\\[1ex]
          \includegraphics[trim=0mm 1mm -1mm 0mm,clip,scale=0.8]{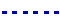}\:$\int_{0}^{L}\sintL\,(\sdam')^2\dx$
        \end{minipage}
        }}
    \end{overpic}
    \caption{$\ell_1$ = 0.05, $\ell_2$ = 0.025}
    \label{fig_DD1} 
  \end{subfigure}
  \hfill
  \begin{subfigure}[b]{0.35\linewidth}
    \begin{overpic}[scale=0.85]{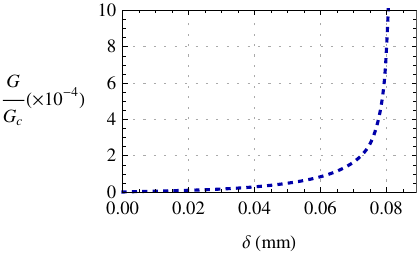}
      \put(40,18){\rotatebox{1}{\normalsize\color{blue}$\boldsymbol\rightarrow$}}
    \end{overpic}
    \caption{}
    \label{fig_DD2} 
  \end{subfigure}
  \caption{(\subref{fig_DD1}) crack density function decomposition between the local and gradient contributions for \eqref{mod_D1} and $\sintL = \SI{1}{mm}$. (\subref{fig_DD2}) zoom on the gradient term. The arrows denote the evolution direction.}
  \label{fig_DD} 
\end{figure}

To further understand the model behaviour, let us highlight in \Fig{fig_DD} the 
energetic contributions to the crack density function during the crack evolution for \eqref{mod_D1}. As noticeable from \Fig{fig_DD1}, the dissipated fracture energy is initially linearly increasing with respect to the crack opening, consistent with the expected Dugdale cohesive response up to the critical crack opening. During this stage, the local term is the dominant contribution. \Fig{fig_DD2} highlights the initial negligible contribution of the gradient term to the crack density function. For a further phase-field evolution up to $\sdammax=1$, the analytical crack-opening computed in \eqref{def_sope} decreases to zero, with a decrease of the local term contribution and a simultaneous significant increase of the gradient term such that the total dissipated fracture energy remains constant.
It is worth remarking that this second stage, although analytically well defined, is not observable for a monotonically increasing bar-end displacement, since, after the critical crack opening, only the occurrence of a brutal crack is compatible through a snap-back response.
It is worth noting that the dissipated fracture energy is always non decreasing, therefore satisfying the second law of thermodynamics, although its single contributions, as the local term, are not. This fact allow us to remark that the gradient phase-field term must be considered as a contribution to the dissipated energy and not to the free energy, as instead done in many works.

%

\begin{figure}[h]
\centering
\begin{subfigure}[b]{\linewidth}
  \centering
  \begin{overpic}[width=0.8\linewidth]{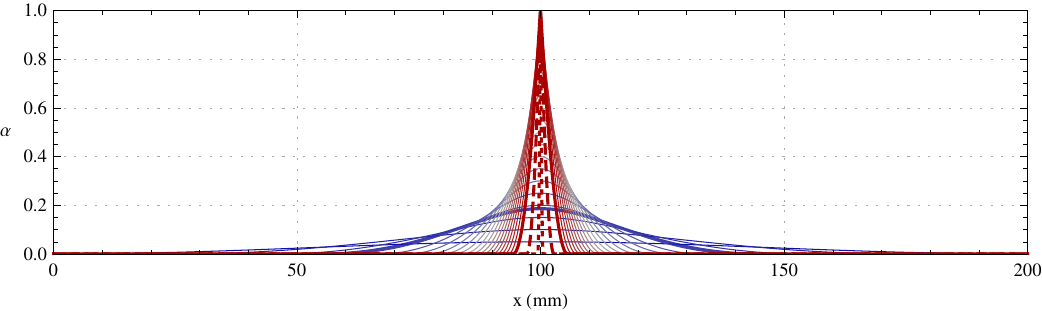}
    \put(80,25){\damBar}
    \put(75,15){\lDILeg}
  \end{overpic}
  \caption{}
  \label{DProfDPlot1} 
\end{subfigure}
\\[2ex]
\begin{subfigure}[b]{\linewidth}
  \centering
  \begin{overpic}[width=0.8\linewidth]{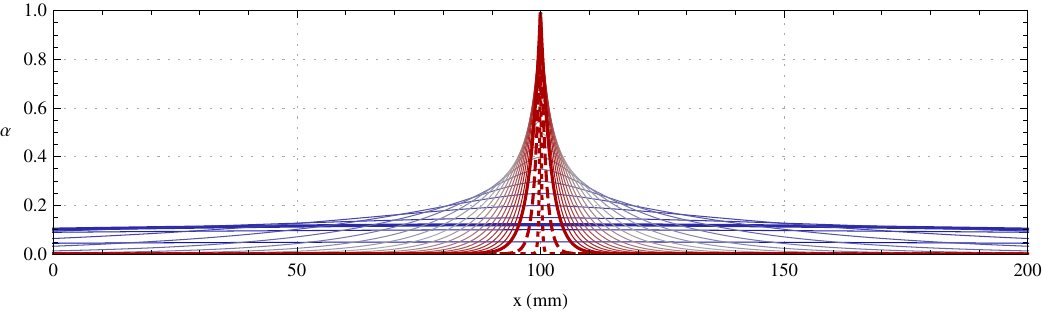}
    \put(80,25){\damBar}
    \put(75,15){\lDILeg}
  \end{overpic}
  \caption{}
  \label{DProfDPlot2} 
\end{subfigure}
\\[2ex]
\begin{subfigure}[b]{\linewidth}
  \centering
  \scriptsize
  \begin{overpic}[width=0.8\linewidth]{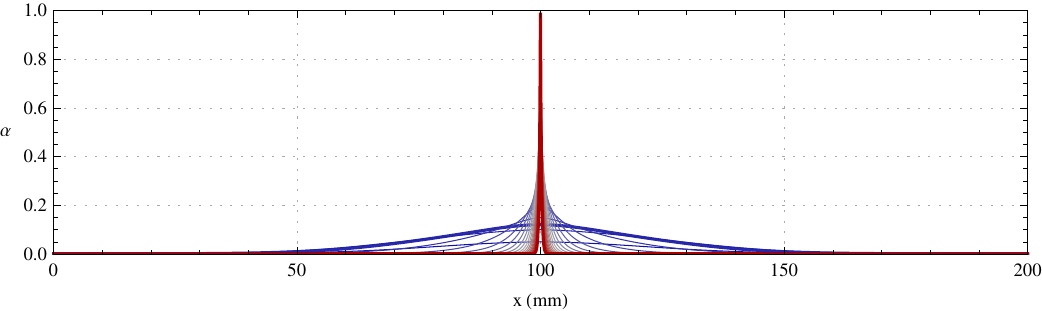}
    \put(80,22){\damBar}
    \put(82,17){$\sintL = \SI{0.1}{mm}$}
  \end{overpic}
  \caption{$\ell=0.1$}
  \label{DProfDPlotScaled2} 
\end{subfigure}
    \caption{Phase-field profiles at given maximum phase-field states for the Dugdale models: (\subref{DProfDPlot1}) for \eqref{mod_D1} with $\sintL=\SI{1}{mm}$, (\subref{DProfDPlot2}) for \eqref{mod_D2} with $\sintL=\SI{1}{mm}$ and (\subref{DProfDPlotScaled2}) for \eqref{mod_D2} with $\sintL=\SI{0.1}{mm}$. In all subfigures the phase-field profile associated with the critical crack opening displacement is highglighted by a thick line.
    In (\subref{DProfDPlot1}) and (\subref{DProfDPlot2}) the ultimate phase-field profiles associated with the three internal lengths~\eqref{def_lD} are also highlighted.}
  \label{fig_DamProfD} 
\end{figure}

\begin{figure}[h]
\centering
\footnotesize
  \begin{overpic}[scale=0.9]{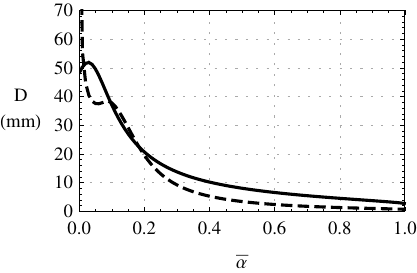}
    \put(22,58){$\uparrow +\infty$}
    \put(70,52){\DMLeg}    
  \end{overpic}
  \caption{Half-width of the support of the phase-field profiles for the two Dugdale models with respect to the fracture evolution ($\sintL = \SI{0.5}{mm}$ for \eqref{mod_D1} and $\sintL = \SI{0.025}{mm}$ for \eqref{mod_D2}).}
  \label{fig_DDevo} 
\end{figure}

\Fig{fig_DamProfD} reports the phase-field profiles evolutions, obtained from \eqref{eq_damprof}, for the two Dugdale models by considering the following internal lengths
\begin{equation} \label{def_lD} 
  \sintL = \cbra{1,0.5,0.1} \quad \text{(mm)} \,.
\end{equation} 
Specifically Fig.s~\ref{DProfDPlot1} and \ref{DProfDPlot2} report the phase-field evolution for the largest internal length together with the ultimate phase-field profiles associated with the smaller internal lengths.
\Fig{DProfDPlotScaled2} highlights the phase-field evolution for the smallest internal length considered in \eqref{def_lD}.
In all three Fig.s \ref{DProfDPlot1}-\ref{DProfDPlotScaled2}, the phase-field profile associated with the critical crack opening displacement is also highlighted by a thick curve.
The trends of the half-width localization support~\eqref{eq_D} with respect to the fracture evolution are highlighted in \Fig{fig_DDevo}. Two different values of the internal length have been considered, one for each model, in order to let the dimensions be comparable. Indeed, \eqref{mod_D1} admits a finite value for the ultimate crack opening whereas \eqref{mod_D2} does not.

It is then immediately clear that although the two Dugdale models describe the same cohesive fracture response and global response, apart from the critical opening displacement values, the behaviour of their corresponding phase-field profiles is not the same. Clearly, for \eqref{mod_D2}, the phase-field at the boundaries has been allowed to assume values different from zero. The phase-field profiles for model \eqref{mod_D1} instead have a finite support.
For both models we notice that the localization support significantly shrinks after the critical phase-field value as $\sdammax$ increases. This also explains the behaviour of the crack density function observed in \Fig{fig_DD1}.

\begin{figure}[h]
\centering
\begin{subfigure}[b]{\linewidth}
  \centering
  \begin{overpic}[width=0.8\linewidth]{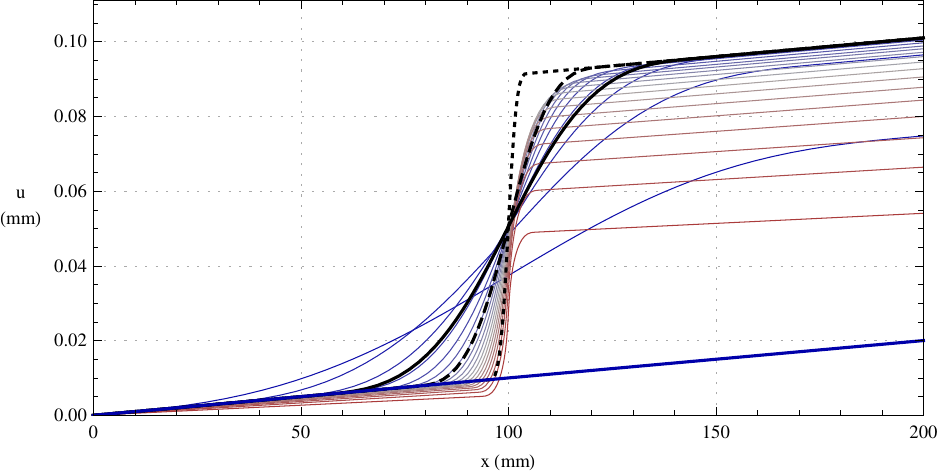}
    \put(11,45){\damBar}
    \put(11,35){\lDILeg[\:\sdam^* = 0.186]}
  \end{overpic}
  \caption{}
  \label{UProfilesDPlot1} 
\end{subfigure}
\\[2ex]
\begin{subfigure}[b]{\linewidth}
  \centering
  \begin{overpic}[width=0.8\linewidth]{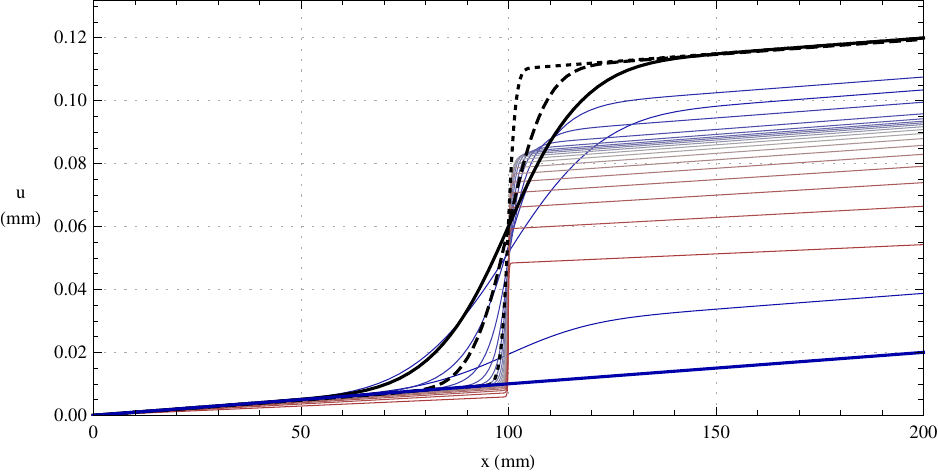}
    \put(11,45){\damBar}
    \put(11,35){\lDILeg[\:\sdam^* = 0.120]}
  \end{overpic}
  \caption{}
  \label{UProfilesDPlot2} 
\end{subfigure}
  \caption{Displacement profiles (\subref{UProfilesDPlot1}) and (\subref{UProfilesDPlot2}) at given maximum phase-field states for, respectively, \eqref{mod_D1} and \eqref{mod_D2} with $\sintL = \SI{1}{mm}$. The critical displacement profiles at $\sdammax = \sdam^* = 0.186$ for \eqref{mod_D1} and $\sdammax = \sdam^* = 0.120$  for \eqref{mod_D2} for smaller internal lengths are also highlighted corresponding to the critical state, where a snap-back of the response occurs. The displacement profiles corresponding to the elastic limit are highlighted with a thick blue line.}
  \label{fig_UProfD} 
\end{figure}

The evolution of the displacement profiles for \eqref{mod_D1} and \eqref{mod_D2} are reported in \Fig{fig_UProfD}. Therein, the convergence of the smooth displacement profiles towards profiles with a discontinuity (the crack) is noticeable as the internal length decreases.

We note that the reconstruction problem for the Dugdale cohesive law was first addressed in \cite[Section 7.1]{Conti2016} and further investigated in \cite{Lammen2025} in a finite strain setting for phase-field models based on the truncated degradation function \eqref{e:f eps tilde}.
For such models, a further regularization of the truncated degradation function is required for numerical implementation and analyses~\cite{Freddi2017,Lammen2025}.
Moreover, the critical stress depends on the internal regularization length, and the global mechanical response characteristic of a Dugdale-type cohesive law is recovered only in the limit as this length tends to zero (see, for instance, Fig.~2 in \cite{Lammen2025}).

Differently, the present work adopts the smooth degradation function \eqref{e:f eps}, which is directly implementable for numerical simulations. As a result, the (sufficiently small) internal regularization length of the phase-field model does not affect the global mechanical response, but only influences the shape of the phase-field and displacement profiles.

%



\FloatBarrier
\section{Conclusions and future perspectives}
\label{sec_conc} 

This work concludes a three-part effort aimed at developing a consistent and unified framework for phase-field modeling of cohesive fracture. The first part of this work \cite{Alessi2025a} established a $\Gamma$-convergence result for a broad class of phase-field energy functionals, encompassing several existing models as \cite{Conti2016,Wu2017,Wu2018b,Feng2021,Lammen2023,Lammen2025}, and demonstrated how an appropriate scaling can generalize the classical Ambrosio-Tortorelli model for brittle fracture. Concerning the models of \cite{Wu2017,Wu2018b,Feng2021}, characterized by a smooth (not-truncated) degradation function of the kind \eqref{e:f eps}, the $\Gamma$-convergence result was expected by numerical investigations but never rigorously proved before this work.
The second part \cite{Alessi2025b} introduced a rigorous analytical methodology, mathematically validating and extending the results in \cite{Feng2021}, for constructing phase-field models tailored to specific cohesive laws, either by prescribing the degradation function and deriving the local phase-field dissipation potential or vice versa. 
This approach enables the generation of novel cohesive phase-field fracture models, including the derivation of multiple models yielding identical cohesive fracture responses but differing in the evolution of their phase-field and displacement profiles.
The capability of generating entire families of phase-field models associated with the same cohesive fracture response deserves undoubtedly to be further investigated, for instance for understanding the possibility of deriving models which satisfy \eqref{def_irr} during the entire evolution.
Building on these theoretical foundations, the present third-part paper examined the mechanical responses of such models under a simple one-dimensional setup. By adopting an engineering-oriented perspective, it provided a practical interpretation of the framework and connected the theoretical insights established in the previous parts to physical applications. 
Possible future developments could include the extension of the present framework to higher dimensional domains, a path already undertaken for instance in \cite{Conti2024,Colasanto2025,Conti2025}, together with the numerical implementation as in \cite{Kumar2023}, the modeling of irreversibility and crack-closure behaviour on the wake of what done in \cite{Crismale2016a,Bonacini2021} for the embedding of fatigue effects~\cite{Alessi2019c,Alessi2023} and the modeling of general mixed-mode cohesive laws \cite{Feng2023}.


\FloatBarrier
\section*{Acknowledgments}

The first author has been supported by the European Union - Next Generation EU, Mission 4 Component 1 CUP G53D23001140006, codice 20229BM9EL, PRIN2022 project: "NutShell - NUmerical modelling and opTimisation of SHELL Structures Against Fracture and Fatigue with Experimental Validations". The first author also acknowledges the Italian National Group of Mathematical Physics INdAM-GNFM. 

The second and third authors have been supported by the European Union - Next Generation EU, Mission 4 Component 1 CUP B53D23009310006, codice 2022J4FYNJ, PRIN2022 project "Variational methods for stationary and evolution problems with singularities and interfaces". The second and third authors are members of
GNAMPA - INdAM.


\bibliographystyle{elsarticle-num-names}
\bibliography{%
  bibliography-mod}

\begin{thebibliography}{55}
\expandafter\ifx\csname natexlab\endcsname\relax\def\natexlab#1{#1}\fi
\providecommand{\url}[1]{\texttt{#1}}
\providecommand{\href}[2]{#2}
\providecommand{\path}[1]{#1}
\providecommand{\DOIprefix}{doi:}
\providecommand{\ArXivprefix}{arXiv:}
\providecommand{\URLprefix}{URL: }
\providecommand{\Pubmedprefix}{pmid:}
\providecommand{\doi}[1]{\href{http://dx.doi.org/#1}{\path{#1}}}
\providecommand{\Pubmed}[1]{\href{pmid:#1}{\path{#1}}}
\providecommand{\bibinfo}[2]{#2}
\ifx\xfnm\relax \def\xfnm[#1]{\unskip,\space#1}\fi
\bibitem[{Griffith(1921)}]{Griffith1921}
\bibinfo{author}{A.~A. Griffith},
\newblock \bibinfo{title}{{The phenomena of rupture and flow in solids}},
\newblock \bibinfo{journal}{Philosophical Transactions of the Royal Society of
  London} \bibinfo{volume}{A221} (\bibinfo{year}{1921})
  \bibinfo{pages}{163--198}. \URLprefix
  \url{http://www.jstor.org/stable/91192}.
\bibitem[{Marigo(2010)}]{marigo2010}
\bibinfo{author}{J.-J. Marigo},
\newblock \bibinfo{title}{{Initiation of Cracks in Griffith's Theory: An
  Argument of Continuity in Favor of Global Minimization}},
\newblock \bibinfo{journal}{Journal of Nonlinear Science} \bibinfo{volume}{20}
  (\bibinfo{year}{2010}) \bibinfo{pages}{831--868}.
  \DOIprefix\doi{10.1007/s00332-010-9074-x}.
\bibitem[{Tann{\'{e}} et~al.(2018)Tann{\'{e}}, Li, Bourdin, Marigo, and
  Maurini}]{Tanne2018}
\bibinfo{author}{E.~Tann{\'{e}}}, \bibinfo{author}{T.~Li},
  \bibinfo{author}{B.~Bourdin}, \bibinfo{author}{J.-J. Marigo},
  \bibinfo{author}{C.~Maurini},
\newblock \bibinfo{title}{{Crack nucleation in variational phase-field models
  of brittle fracture}},
\newblock \bibinfo{journal}{Journal of the Mechanics and Physics of Solids}
  \bibinfo{volume}{110} (\bibinfo{year}{2018}) \bibinfo{pages}{80--99}.
\bibitem[{Kumar et~al.(2020)Kumar, Bourdin, Francfort, and
  Lopez-Pamies}]{Kumar2020}
\bibinfo{author}{A.~Kumar}, \bibinfo{author}{B.~Bourdin},
  \bibinfo{author}{G.~A. Francfort}, \bibinfo{author}{O.~Lopez-Pamies},
\newblock \bibinfo{title}{{Revisiting nucleation in the phase-Field approach to
  brittle fracture}},
\newblock \bibinfo{journal}{Journal of the Mechanics and Physics of Solids}
  \bibinfo{volume}{142} (\bibinfo{year}{2020}) \bibinfo{pages}{104027}.
  \URLprefix \url{https://doi.org/10.1016/j.jmps.2020.104027}.
  \DOIprefix\doi{10.1016/j.jmps.2020.104027}.
\bibitem[{Marigo(2023)}]{Marigo2023}
\bibinfo{author}{J.-J. Marigo},
\newblock \bibinfo{title}{{Modelling of fracture by cohesive force models: A
  path to pursue}},
\newblock \bibinfo{journal}{European Journal of Mechanics, A/Solids}
  \bibinfo{volume}{102} (\bibinfo{year}{2023}).
  \DOIprefix\doi{10.1016/j.euromechsol.2023.105088}.
\bibitem[{Dugdale(1960)}]{Dugdale1960}
\bibinfo{author}{D.~Dugdale},
\newblock \bibinfo{title}{{Yielding of steel sheets containing slits}},
\newblock \bibinfo{journal}{Journal of Mechanics Physics of Solids}
  \bibinfo{volume}{8} (\bibinfo{year}{1960}) \bibinfo{pages}{100--104}.
  \DOIprefix\doi{10.1016/0022-5096(60)90013-2}.
\bibitem[{Barenblatt(1962)}]{Barenblatt1962}
\bibinfo{author}{G.~I. Barenblatt},
\newblock \bibinfo{title}{{The mathematical theory of equilibrium of cracks in
  brittle fracture}},
\newblock \bibinfo{journal}{Advances in Applied Mechanics} \bibinfo{volume}{7}
  (\bibinfo{year}{1962}) \bibinfo{pages}{55--129}.
\bibitem[{Francfort and Marigo(1998)}]{Marigo1998}
\bibinfo{author}{G.~A. Francfort}, \bibinfo{author}{J.-J. Marigo},
\newblock \bibinfo{title}{{Revisiting brittle fracture as an energy
  minimization problem}},
\newblock \bibinfo{journal}{Journal of the Mechanics and Physics of Solids}
  \bibinfo{volume}{46} (\bibinfo{year}{1998}) \bibinfo{pages}{1319--1342}.
  \URLprefix
  \url{http://www.sciencedirect.com/science/article/B6TXB-3XDH3V4-2/2/6078bc5b38861a78a6b1044a60e0a4bb}.
  \DOIprefix\doi{10.1016/S0022-5096(98)00034-9}.
\bibitem[{Focardi(2001)}]{Focardi2001}
\bibinfo{author}{M.~Focardi},
\newblock \bibinfo{title}{{On the variational approximation of
  free-discontinuity problems in the vectorial case}},
\newblock \bibinfo{journal}{Mathematical Models and Methods in Applied
  Sciences} \bibinfo{volume}{11} (\bibinfo{year}{2001})
  \bibinfo{pages}{663--684}. \URLprefix
  \url{http://dx.doi.org/10.1142/S0218202501001045
  https://www.worldscientific.com/doi/10.1142/S0218202501001045}.
  \DOIprefix\doi{10.1142/S0218202501001045}.
\bibitem[{Bourdin et~al.(2008)Bourdin, Francfort, and Marigo}]{Bourdin2008}
\bibinfo{author}{B.~Bourdin}, \bibinfo{author}{G.~A. Francfort},
  \bibinfo{author}{J.-J. Marigo},
\newblock \bibinfo{title}{{The Variational Approach to Fracture}},
\newblock \bibinfo{journal}{Journal of Elasticity} \bibinfo{volume}{91}
  (\bibinfo{year}{2008}) \bibinfo{pages}{5--148}. \URLprefix
  \url{http://www.springerlink.com/content/31370771n83r2637}.
  \DOIprefix\doi{10.1007/s10659-007-9107-3}.
\bibitem[{Bourdin et~al.(2000)Bourdin, Francfort, and Marigo}]{Bourdin2000b}
\bibinfo{author}{B.~Bourdin}, \bibinfo{author}{G.~A. Francfort},
  \bibinfo{author}{J.-J. Marigo},
\newblock \bibinfo{title}{{Numerical experiments in revisited brittle
  fracture}},
\newblock \bibinfo{journal}{Journal of the Mechanics and Physics of Solids}
  \bibinfo{volume}{48} (\bibinfo{year}{2000}) \bibinfo{pages}{797--826}.
  \URLprefix
  \url{http://linkinghub.elsevier.com/retrieve/pii/S0022509699000289}.
  \DOIprefix\doi{10.1016/S0022-5096(99)00028-9}.
\bibitem[{Ambrosio and Tortorelli(1992)}]{Ambrosio1992}
\bibinfo{author}{L.~Ambrosio}, \bibinfo{author}{V.~M. Tortorelli},
\newblock \bibinfo{title}{{On the approximation of free discontinuity
  problems}},
\newblock \bibinfo{journal}{Boll. Un. Mat. Ital.} \bibinfo{volume}{6-B}
  (\bibinfo{year}{1992}) \bibinfo{pages}{105--123}.
\bibitem[{Braides(1998)}]{Braides1998}
\bibinfo{author}{A.~Braides}, \bibinfo{title}{{Approximation of
  Free-Discontinuity Problems}}, volume \bibinfo{volume}{1694} of
  \textit{\bibinfo{series}{Lecture Notes in Mathematics}},
  \bibinfo{publisher}{Springer Berlin Heidelberg}, \bibinfo{address}{Berlin,
  Heidelberg}, \bibinfo{year}{1998}. \URLprefix
  \url{http://link.springer.com/10.1007/BFb0097344}.
  \DOIprefix\doi{10.1007/BFb0097344}.
\bibitem[{Bourdin et~al.(2014)Bourdin, Marigo, Maurini, and
  Sicsic}]{Bourdin2014}
\bibinfo{author}{B.~Bourdin}, \bibinfo{author}{J.-J. Marigo},
  \bibinfo{author}{C.~Maurini}, \bibinfo{author}{P.~Sicsic},
\newblock \bibinfo{title}{{Morphogenesis and Propagation of Complex Cracks
  Induced by Thermal Shocks}},
\newblock \bibinfo{journal}{Physical Review Letters} \bibinfo{volume}{112}
  (\bibinfo{year}{2014}) \bibinfo{pages}{014301}. \URLprefix
  \url{http://link.aps.org/doi/10.1103/PhysRevLett.112.014301}.
  \DOIprefix\doi{10.1103/PhysRevLett.112.014301}.
\bibitem[{Mesgarnejad et~al.(2015)Mesgarnejad, Bourdin, and
  Khonsari}]{Mesgarnejad2015}
\bibinfo{author}{A.~Mesgarnejad}, \bibinfo{author}{B.~Bourdin},
  \bibinfo{author}{M.~M. Khonsari},
\newblock \bibinfo{title}{{Validation simulations for the variational approach
  to fracture}},
\newblock \bibinfo{journal}{Computer Methods in Applied Mechanics and
  Engineering} \bibinfo{volume}{290} (\bibinfo{year}{2015})
  \bibinfo{pages}{420--437}. \URLprefix
  \url{http://www.sciencedirect.com/science/article/pii/S004578251400423X}.
  \DOIprefix\doi{10.1016/j.cma.2014.10.052}.
\bibitem[{Lemaitre and Chaboche(1985)}]{LemaitreB1985}
\bibinfo{author}{J.~Lemaitre}, \bibinfo{author}{J.-L. Chaboche},
  \bibinfo{title}{{Mechanics of Solid Materials}},
  \bibinfo{publisher}{Cambridge University Press}, \bibinfo{year}{1985}.
  \URLprefix \url{http://books.google.it/books?id=7PCncQAACAAJ}.
\bibitem[{Pham et~al.(2011)Pham, Amor, Marigo, and Maurini}]{Pham2011}
\bibinfo{author}{K.~Pham}, \bibinfo{author}{H.~Amor}, \bibinfo{author}{J.-J.
  Marigo}, \bibinfo{author}{C.~Maurini},
\newblock \bibinfo{title}{{Gradient Damage Models and Their Use to Approximate
  Brittle Fracture}},
\newblock \bibinfo{journal}{International Journal of Damage Mechanics}
  \bibinfo{volume}{20} (\bibinfo{year}{2011}) \bibinfo{pages}{618--652}.
  \DOIprefix\doi{10.1177/1056789510386852}.
\bibitem[{Pham and Marigo(2010{\natexlab{a}})}]{Pham2010c}
\bibinfo{author}{K.~Pham}, \bibinfo{author}{J.-J. Marigo},
\newblock \bibinfo{title}{{The variational approach to damage: II. The gradient
  damage models}},
\newblock \bibinfo{journal}{Comptes Rendus Mecanique} \bibinfo{volume}{338}
  (\bibinfo{year}{2010}{\natexlab{a}}) \bibinfo{pages}{199--206}.
  \DOIprefix\doi{10.1016/j.crme.2010.03.012}.
\bibitem[{Pham and Marigo(2010{\natexlab{b}})}]{Pham2010a}
\bibinfo{author}{K.~Pham}, \bibinfo{author}{J.-J. Marigo},
\newblock \bibinfo{title}{{The variational approach to damage: I. The
  foundations}},
\newblock \bibinfo{journal}{Comptes Rendus Mecanique} \bibinfo{volume}{338}
  (\bibinfo{year}{2010}{\natexlab{b}}) \bibinfo{pages}{191--198}.
  \DOIprefix\doi{10.1016/j.crme.2010.03.009}.
\bibitem[{Miehe et~al.(2010)Miehe, Welschinger, and Hofacker}]{Miehe2010a}
\bibinfo{author}{C.~Miehe}, \bibinfo{author}{F.~R. Welschinger},
  \bibinfo{author}{M.~Hofacker},
\newblock \bibinfo{title}{{Thermodynamically consistent phase-field models of
  fracture: Variational principles and multi-field FE implementations}},
\newblock \bibinfo{journal}{International Journal for Numerical Methods in
  Engineering} \bibinfo{volume}{83} (\bibinfo{year}{2010})
  \bibinfo{pages}{1273--1311}. \URLprefix
  \url{http://doi.wiley.com/10.1002/nme.2861}.
  \DOIprefix\doi{10.1002/nme.2861}.
\bibitem[{Marigo et~al.(2016)Marigo, Maurini, and Pham}]{Marigo2016}
\bibinfo{author}{J.-J. Marigo}, \bibinfo{author}{C.~Maurini},
  \bibinfo{author}{K.~Pham},
\newblock \bibinfo{title}{{An overview of the modelling of fracture by gradient
  damage models}},
\newblock \bibinfo{journal}{Meccanica} \bibinfo{volume}{51}
  (\bibinfo{year}{2016}) \bibinfo{pages}{3107--3128}. \URLprefix
  \url{http://link.springer.com/10.1007/s11012-016-0538-4
  https://hal-ensmp.archives-ouvertes.fr/hal-01321465/}.
  \DOIprefix\doi{10.1177/1056789510386852}.
\bibitem[{Lopez-Pamies et~al.(2024)Lopez-Pamies, Dolbow, Francfort, and
  Larsen}]{Lopez-Pamies2024}
\bibinfo{author}{O.~Lopez-Pamies}, \bibinfo{author}{J.~E. Dolbow},
  \bibinfo{author}{G.~A. Francfort}, \bibinfo{author}{C.~J. Larsen},
\newblock \bibinfo{title}{{Classical variational phase-field models cannot
  predict fracture nucleation}}  (\bibinfo{year}{2024}). \URLprefix
  \url{http://arxiv.org/abs/2409.00242}.
  \href{http://arxiv.org/abs/2409.00242}{{\tt arXiv:2409.00242}}.
\bibitem[{Lorentz et~al.(2011)Lorentz, Cuvilliez, and
  Kazymyrenko}]{Lorentz2011}
\bibinfo{author}{E.~Lorentz}, \bibinfo{author}{S.~Cuvilliez},
  \bibinfo{author}{K.~Kazymyrenko},
\newblock \bibinfo{title}{{Convergence of a gradient damage model toward a
  cohesive zone model}},
\newblock \bibinfo{journal}{Comptes Rendus Mecanique} \bibinfo{volume}{339}
  (\bibinfo{year}{2011}) \bibinfo{pages}{20--26}. \URLprefix
  \url{http://www.sciencedirect.com/science/article/pii/S1631072110001671}.
  \DOIprefix\doi{10.1016/j.crme.2010.10.010}.
\bibitem[{Lorentz and Godard(2011)}]{Lorentz2011a}
\bibinfo{author}{E.~Lorentz}, \bibinfo{author}{V.~Godard},
\newblock \bibinfo{title}{{Gradient damage models: Toward full-scale
  computations}},
\newblock \bibinfo{journal}{Computer Methods in Applied Mechanics and
  Engineering} \bibinfo{volume}{200} (\bibinfo{year}{2011})
  \bibinfo{pages}{1927--1944}. \URLprefix
  \url{http://www.sciencedirect.com/science/article/pii/S0045782510001921}.
  \DOIprefix\doi{10.1016/j.cma.2010.06.025}.
\bibitem[{Conti et~al.(2016)Conti, Focardi, and Iurlano}]{Conti2016}
\bibinfo{author}{S.~Conti}, \bibinfo{author}{M.~Focardi},
  \bibinfo{author}{F.~Iurlano},
\newblock \bibinfo{title}{{Phase field approximation of cohesive fracture
  models}},
\newblock \bibinfo{journal}{Annales de l'Institut Henri Poincar{\'{e}} C,
  Analyse non lin{\'{e}}aire} \bibinfo{volume}{33} (\bibinfo{year}{2016})
  \bibinfo{pages}{1033--1067}. \URLprefix
  \url{http://www.sciencedirect.com/science/article/pii/S0294144915000360
  https://ems.press/doi/10.1016/j.anihpc.2015.02.001}.
  \DOIprefix\doi{10.1016/j.anihpc.2015.02.001}.
  \href{http://arxiv.org/abs/arXiv:1405.6883v1}{{\tt arXiv:arXiv:1405.6883v1}}.
\bibitem[{Freddi and Iurlano(2017)}]{Freddi2017}
\bibinfo{author}{F.~Freddi}, \bibinfo{author}{F.~Iurlano},
\newblock \bibinfo{title}{{Numerical insight of a variational smeared approach
  to cohesive fracture}},
\newblock \bibinfo{journal}{Journal of the Mechanics and Physics of Solids}
  \bibinfo{volume}{98} (\bibinfo{year}{2017}) \bibinfo{pages}{156--171}.
  \URLprefix
  \url{http://www.sciencedirect.com/science/article/pii/S0022509616304860}.
  \DOIprefix\doi{10.1016/j.jmps.2016.09.003}.
\bibitem[{Wu(2017)}]{Wu2017}
\bibinfo{author}{J.-Y. Wu},
\newblock \bibinfo{title}{{A unified phase-field theory for the mechanics of
  damage and quasi-brittle failure}},
\newblock \bibinfo{journal}{Journal of the Mechanics and Physics of Solids}
  \bibinfo{volume}{103} (\bibinfo{year}{2017}) \bibinfo{pages}{72--99}.
  \URLprefix \url{http://dx.doi.org/10.1016/j.jmps.2017.03.015
  http://www.sciencedirect.com/science/article/pii/S0022509616308341}.
  \DOIprefix\doi{10.1016/j.jmps.2017.03.015}.
\bibitem[{Wu and Nguyen(2018)}]{Wu2018b}
\bibinfo{author}{J.-Y. Wu}, \bibinfo{author}{V.~P. Nguyen},
\newblock \bibinfo{title}{{A length scale insensitive phase-field damage model
  for brittle fracture}},
\newblock \bibinfo{journal}{Journal of the Mechanics and Physics of Solids}
  \bibinfo{volume}{119} (\bibinfo{year}{2018}) \bibinfo{pages}{20--42}.
  \URLprefix \url{https://doi.org/10.1016/j.jmps.2018.06.006}.
  \DOIprefix\doi{10.1016/j.jmps.2018.06.006}.
\bibitem[{Feng et~al.(2021)Feng, Fan, and Li}]{Feng2021}
\bibinfo{author}{Y.~Feng}, \bibinfo{author}{J.~Fan}, \bibinfo{author}{J.~Li},
\newblock \bibinfo{title}{{Endowing explicit cohesive laws to the phase-field
  fracture theory}},
\newblock \bibinfo{journal}{Journal of the Mechanics and Physics of Solids}
  \bibinfo{volume}{152} (\bibinfo{year}{2021}) \bibinfo{pages}{104464}.
  \URLprefix \url{https://doi.org/10.1016/j.jmps.2021.104464}.
  \DOIprefix\doi{10.1016/j.jmps.2021.104464}.
\bibitem[{Lammen et~al.(2023)Lammen, Conti, and Mosler}]{Lammen2023}
\bibinfo{author}{H.~Lammen}, \bibinfo{author}{S.~Conti},
  \bibinfo{author}{J.~Mosler},
\newblock \bibinfo{title}{{A finite deformation phase field model suitable for
  cohesive fracture}},
\newblock \bibinfo{journal}{Journal of the Mechanics and Physics of Solids}
  \bibinfo{volume}{178} (\bibinfo{year}{2023}) \bibinfo{pages}{105349}.
  \URLprefix
  \url{https://linkinghub.elsevier.com/retrieve/pii/S0022509623001539}.
  \DOIprefix\doi{10.1016/j.jmps.2023.105349}.
\bibitem[{Lammen et~al.(2025)Lammen, Conti, and Mosler}]{Lammen2025}
\bibinfo{author}{H.~Lammen}, \bibinfo{author}{S.~Conti},
  \bibinfo{author}{J.~Mosler},
\newblock \bibinfo{title}{{Approximating arbitrary traction–separation-laws
  by means of phase-field theory — Mathematical foundation and numerical
  implementation}},
\newblock \bibinfo{journal}{Journal of the Mechanics and Physics of Solids}
  \bibinfo{volume}{197} (\bibinfo{year}{2025}) \bibinfo{pages}{106038}.
  \URLprefix
  \url{https://linkinghub.elsevier.com/retrieve/pii/S0022509625000146}.
  \DOIprefix\doi{10.1016/j.jmps.2025.106038}.
\bibitem[{Alessi et~al.(2025{\natexlab{a}})Alessi, Colasanto, and
  Focardi}]{Alessi2025a}
\bibinfo{author}{R.~Alessi}, \bibinfo{author}{F.~Colasanto},
  \bibinfo{author}{M.~Focardi},
\newblock \bibinfo{title}{{Phase-field modelling of cohesive fracture. Part I:
  $\Gamma$-convergence results}},
\newblock \bibinfo{journal}{~\!\!} \bibinfo{volume}{submitted}
  (\bibinfo{year}{2025}{\natexlab{a}}).
\bibitem[{Alessi et~al.(2025{\natexlab{b}})Alessi, Colasanto, and
  Focardi}]{Alessi2025b}
\bibinfo{author}{R.~Alessi}, \bibinfo{author}{F.~Colasanto},
  \bibinfo{author}{M.~Focardi},
\newblock \bibinfo{title}{{Phase-field modelling of cohesive fracture. Part II:
  Reconstruction of the cohesive law}},
\newblock \bibinfo{journal}{~\!\!} \bibinfo{volume}{submitted}
  (\bibinfo{year}{2025}{\natexlab{b}}).
\bibitem[{Alessi et~al.(2014)Alessi, Marigo, and Vidoli}]{Alessi2014}
\bibinfo{author}{R.~Alessi}, \bibinfo{author}{J.-J. Marigo},
  \bibinfo{author}{S.~Vidoli},
\newblock \bibinfo{title}{{Gradient Damage Models Coupled with Plasticity and
  Nucleation of Cohesive Cracks}},
\newblock \bibinfo{journal}{Archive for Rational Mechanics and Analysis}
  \bibinfo{volume}{214} (\bibinfo{year}{2014}) \bibinfo{pages}{575--615}.
  \URLprefix \url{http://link.springer.com/10.1007/s00205-014-0763-8}.
  \DOIprefix\doi{10.1007/s00205-014-0763-8}.
\bibitem[{Alessi et~al.(2015)Alessi, Marigo, and Vidoli}]{Alessi2015}
\bibinfo{author}{R.~Alessi}, \bibinfo{author}{J.-J. Marigo},
  \bibinfo{author}{S.~Vidoli},
\newblock \bibinfo{title}{{Gradient damage models coupled with plasticity:
  Variational formulation and main properties}},
\newblock \bibinfo{journal}{Mechanics of Materials} \bibinfo{volume}{80}
  (\bibinfo{year}{2015}) \bibinfo{pages}{351--367}. \URLprefix
  \url{http://www.sciencedirect.com/science/article/pii/S0167663614000039}.
  \DOIprefix\doi{10.1016/j.mechmat.2013.12.005}.
\bibitem[{Alessi et~al.(2018)Alessi, Marigo, Maurini, and Vidoli}]{Alessi2018a}
\bibinfo{author}{R.~Alessi}, \bibinfo{author}{J.-J. Marigo},
  \bibinfo{author}{C.~Maurini}, \bibinfo{author}{S.~Vidoli},
\newblock \bibinfo{title}{{Coupling damage and plasticity for a phase-field
  regularisation of brittle, cohesive and ductile fracture: One-dimensional
  examples}},
\newblock \bibinfo{journal}{International Journal of Mechanical Sciences}
  \bibinfo{volume}{149} (\bibinfo{year}{2018}) \bibinfo{pages}{559--576}.
  \DOIprefix\doi{10.1016/j.ijmecsci.2017.05.047}.
\bibitem[{Conti et~al.(2024)Conti, Focardi, and Iurlano}]{Conti2024}
\bibinfo{author}{S.~Conti}, \bibinfo{author}{M.~Focardi},
  \bibinfo{author}{F.~Iurlano},
\newblock \bibinfo{title}{{Phase-Field Approximation of a Vectorial,
  Geometrically Nonlinear Cohesive Fracture Energy}},
\newblock \bibinfo{journal}{Archive for Rational Mechanics and Analysis}
  \bibinfo{volume}{248} (\bibinfo{year}{2024}) \bibinfo{pages}{21}. \URLprefix
  \url{https://link.springer.com/10.1007/s00205-024-01962-4}.
  \DOIprefix\doi{10.1007/s00205-024-01962-4}.
\bibitem[{Colasanto(2025)}]{Colasanto2025}
\bibinfo{author}{F.~Colasanto},
\newblock \bibinfo{title}{{Phase field approximation of cohesive functionals in
  the vectorial case}},
\newblock \bibinfo{journal}{In preparation}  (\bibinfo{year}{2025}).
\bibitem[{Conti et~al.(2025)Conti, Focardi, and Iurlano}]{Conti2025}
\bibinfo{author}{S.~Conti}, \bibinfo{author}{M.~Focardi},
  \bibinfo{author}{F.~Iurlano},
\newblock \bibinfo{title}{{Superlinear free-discontinuity models: relaxation
  and phase field approximation}},
\newblock \bibinfo{journal}{In preparation}  (\bibinfo{year}{2025}). \URLprefix
  \url{http://arxiv.org/abs/2505.00852}. \DOIprefix\doi{arXiv.2505.00852}.
  \href{http://arxiv.org/abs/2505.00852}{{\tt arXiv:2505.00852}}.
\bibitem[{Pham and Marigo(2011)}]{Pham2011c}
\bibinfo{author}{K.~Pham}, \bibinfo{author}{J.-J. Marigo},
\newblock \bibinfo{title}{{From the onset of damage to rupture: construction of
  responses with damage localization for a general class of gradient damage
  models}},
\newblock \bibinfo{journal}{Continuum Mechanics and Thermodynamics}
  \bibinfo{volume}{25} (\bibinfo{year}{2011}) \bibinfo{pages}{147--171}.
  \URLprefix \url{http://link.springer.com/10.1007/s00161-011-0228-3}.
  \DOIprefix\doi{10.1007/s00161-011-0228-3}.
\bibitem[{Chen and de~Borst(2021)}]{Chen2021a}
\bibinfo{author}{L.~Chen}, \bibinfo{author}{R.~de~Borst},
\newblock \bibinfo{title}{{Phase-field modelling of cohesive fracture}},
\newblock \bibinfo{journal}{European Journal of Mechanics, A/Solids}
  \bibinfo{volume}{90} (\bibinfo{year}{2021}) \bibinfo{pages}{104343}.
  \URLprefix \url{https://doi.org/10.1016/j.euromechsol.2021.104343}.
  \DOIprefix\doi{10.1016/j.euromechsol.2021.104343}.
\bibitem[{Mielke(2006)}]{Mielke2006}
\bibinfo{author}{A.~Mielke},
\newblock \bibinfo{title}{{A Mathematical Framework for Generalized Standard
  Materials in the Rate-Independent Case}},
\newblock in: \bibinfo{editor}{R.~Helmig}, \bibinfo{editor}{A.~Mielke},
  \bibinfo{editor}{B.~Wohlmuth} (Eds.), \bibinfo{booktitle}{Multifield Problems
  in Solid and Fluid Mechanics}, volume~\bibinfo{volume}{28} of
  \textit{\bibinfo{series}{Lecture Notes in Applied and Computational
  Mechanics}}, \bibinfo{publisher}{Springer Berlin / Heidelberg},
  \bibinfo{year}{2006}, pp. \bibinfo{pages}{399--428}. \URLprefix
  \url{http://dx.doi.org/10.1007/978-3-540-34961-7_12}.
\bibitem[{Mielke and Roub{\'{i}}{\v{c}}ek(2015)}]{Mielke2015}
\bibinfo{author}{A.~Mielke}, \bibinfo{author}{T.~Roub{\'{i}}{\v{c}}ek},
  \bibinfo{title}{{Rate-Independent Systems: Theory and Application}},
  \bibinfo{publisher}{Springer}, \bibinfo{year}{2015}. \URLprefix
  \url{https://books.google.com/books?id=RN3HCgAAQBAJ&pgis=1}.
\bibitem[{Wu(2024)}]{Wu2024a}
\bibinfo{author}{J.-y. Wu},
\newblock \bibinfo{title}{{A generalized phase-field cohesive zone model
  ($\mu$PF-CZM) for fracture}},
\newblock \bibinfo{journal}{Journal of the Mechanics and Physics of Solids}
  \bibinfo{volume}{192} (\bibinfo{year}{2024}).
  \DOIprefix\doi{10.1016/j.jmps.2024.105841}.
  \href{http://arxiv.org/abs/2408.00015}{{\tt arXiv:2408.00015}}.
\bibitem[{Crismale et~al.(2018)Crismale, Lazzaroni, and
  Orlando}]{Crismale2016a}
\bibinfo{author}{V.~Crismale}, \bibinfo{author}{G.~Lazzaroni},
  \bibinfo{author}{G.~Orlando},
\newblock \bibinfo{title}{{Cohesive fracture with irreversibility: quasistatic
  evolution for a model subject to fatigue}},
\newblock \bibinfo{journal}{Mathematical Models and Methods in Applied
  Sciences} \bibinfo{volume}{28} (\bibinfo{year}{2018})
  \bibinfo{pages}{1371--1412}. \URLprefix
  \url{http://cvgmt.sns.it/media/doc/paper/3111/Cri-Laz-Orl-preprint2.pdf}.
  \DOIprefix\doi{10.1142/S0218202518500379}.
\bibitem[{Bonacini et~al.(2021)Bonacini, Conti, and Iurlano}]{Bonacini2021}
\bibinfo{author}{M.~Bonacini}, \bibinfo{author}{S.~Conti},
  \bibinfo{author}{F.~Iurlano},
\newblock \bibinfo{title}{{Cohesive Fracture in 1D: Quasi-static Evolution and
  Derivation from Static Phase-Field Models}},
\newblock \bibinfo{journal}{Archive for Rational Mechanics and Analysis}
  \bibinfo{volume}{239} (\bibinfo{year}{2021}) \bibinfo{pages}{1501--1576}.
  \URLprefix \url{http://link.springer.com/10.1007/s00205-020-01597-1}.
  \DOIprefix\doi{10.1007/s00205-020-01597-1}.
\bibitem[{Guinea et~al.(1994)Guinea, Planas, and Elices}]{Guinea1994}
\bibinfo{author}{G.~V. Guinea}, \bibinfo{author}{J.~Planas},
  \bibinfo{author}{M.~Elices},
\newblock \bibinfo{title}{{A general bilinear fit for the softening curve of
  concrete}},
\newblock \bibinfo{journal}{Materials and Structures} \bibinfo{volume}{27}
  (\bibinfo{year}{1994}) \bibinfo{pages}{99--105}. \URLprefix
  \url{http://link.springer.com/10.1007/BF02472827}.
  \DOIprefix\doi{10.1007/BF02472827}.
\bibitem[{Roesler et~al.(2007)Roesler, Paulino, Park, and
  Gaedicke}]{Roesler2007a}
\bibinfo{author}{J.~Roesler}, \bibinfo{author}{G.~H. Paulino},
  \bibinfo{author}{K.~Park}, \bibinfo{author}{C.~Gaedicke},
\newblock \bibinfo{title}{{Concrete fracture prediction using bilinear
  softening}},
\newblock \bibinfo{journal}{Cement and Concrete Composites}
  \bibinfo{volume}{29} (\bibinfo{year}{2007}) \bibinfo{pages}{300--312}.
  \URLprefix
  \url{https://linkinghub.elsevier.com/retrieve/pii/S0958946506002137}.
  \DOIprefix\doi{10.1016/j.cemconcomp.2006.12.002}.
\bibitem[{Park and Paulino(2011)}]{Park2011}
\bibinfo{author}{K.~Park}, \bibinfo{author}{G.~H. Paulino},
\newblock \bibinfo{title}{{Cohesive Zone Models: A Critical Review of
  Traction-Separation Relationships Across Fracture Surfaces}},
\newblock \bibinfo{journal}{Applied Mechanics Reviews} \bibinfo{volume}{64}
  (\bibinfo{year}{2011}). \URLprefix
  \url{https://asmedigitalcollection.asme.org/appliedmechanicsreviews/article/doi/10.1115/1.4023110/370063/Cohesive-Zone-Models-A-Critical-Review-of}.
  \DOIprefix\doi{10.1115/1.4023110}.
\bibitem[{Dourado et~al.(2012)Dourado, de~Moura, de~Morais, and
  Pereira}]{Dourado2012}
\bibinfo{author}{N.~Dourado}, \bibinfo{author}{M.~de~Moura},
  \bibinfo{author}{A.~de~Morais}, \bibinfo{author}{A.~Pereira},
\newblock \bibinfo{title}{{Bilinear approximations to the mode II delamination
  cohesive law using an inverse method}},
\newblock \bibinfo{journal}{Mechanics of Materials} \bibinfo{volume}{49}
  (\bibinfo{year}{2012}) \bibinfo{pages}{42--50}. \URLprefix
  \url{https://linkinghub.elsevier.com/retrieve/pii/S0167663612000361}.
  \DOIprefix\doi{10.1016/j.mechmat.2012.02.004}.
\bibitem[{Ba{\v{z}}ant et~al.(2002)Ba{\v{z}}ant, Yu, and Zi}]{Bazant2002}
\bibinfo{author}{Z.~P. Ba{\v{z}}ant}, \bibinfo{author}{Q.~Yu},
  \bibinfo{author}{G.~Zi},
\newblock \bibinfo{title}{{Choice of standard fracture test for concrete and
  its statistical evaluation}},
\newblock \bibinfo{journal}{International Journal of Fracture}
  \bibinfo{volume}{118} (\bibinfo{year}{2002}) \bibinfo{pages}{303--337}.
  \URLprefix \url{https://doi.org/10.1023/A:1023399125413}.
  \DOIprefix\doi{10.1023/A:1023399125413}.
\bibitem[{Kumar et~al.(2023)Kumar, Alessi, and Salvati}]{Kumar2023}
\bibinfo{author}{M.~Kumar}, \bibinfo{author}{R.~Alessi},
  \bibinfo{author}{E.~Salvati},
\newblock \bibinfo{title}{{GPFniCS: A generalised phase field method to model
  fracture}},
\newblock \bibinfo{journal}{SoftwareX} \bibinfo{volume}{24}
  (\bibinfo{year}{2023}) \bibinfo{pages}{101594}. \URLprefix
  \url{https://linkinghub.elsevier.com/retrieve/pii/S235271102300290X}.
  \DOIprefix\doi{10.1016/j.softx.2023.101594}.
\bibitem[{Alessi et~al.(2019)Alessi, Crismale, and Orlando}]{Alessi2019c}
\bibinfo{author}{R.~Alessi}, \bibinfo{author}{V.~Crismale},
  \bibinfo{author}{G.~Orlando},
\newblock \bibinfo{title}{{Fatigue Effects in Elastic Materials with
  Variational Damage Models: A Vanishing Viscosity Approach}},
\newblock \bibinfo{journal}{Journal of Nonlinear Science} \bibinfo{volume}{29}
  (\bibinfo{year}{2019}). \DOIprefix\doi{10.1007/s00332-018-9511-9}.
\bibitem[{Alessi and Ulloa(2023)}]{Alessi2023}
\bibinfo{author}{R.~Alessi}, \bibinfo{author}{J.~Ulloa},
\newblock \bibinfo{title}{{Endowing Griffith's fracture theory with the ability
  to describe fatigue cracks}},
\newblock \bibinfo{journal}{Engineering Fracture Mechanics}
  \bibinfo{volume}{281} (\bibinfo{year}{2023}) \bibinfo{pages}{109048}.
  \URLprefix
  \url{https://linkinghub.elsevier.com/retrieve/pii/S0013794423000061}.
  \DOIprefix\doi{10.1016/j.engfracmech.2023.109048}.
\bibitem[{Feng and Li(2023)}]{Feng2023}
\bibinfo{author}{Y.~Feng}, \bibinfo{author}{J.~Li},
\newblock \bibinfo{title}{{A unified regularized variational cohesive fracture
  theory with directional energy decomposition}},
\newblock \bibinfo{journal}{International Journal of Engineering Science}
  \bibinfo{volume}{182} (\bibinfo{year}{2023}) \bibinfo{pages}{103773}.
  \URLprefix \url{https://doi.org/10.1016/j.ijengsci.2022.103773}.
  \DOIprefix\doi{10.1016/j.ijengsci.2022.103773}.

\end{thebibliography}

\end{document}